# The Göttsche conjecture and the Göttsche–Yau–Zaslow formula

## by Francesco Bastianelli and Thomas Dedieu

**Abstract.** The aim of these notes is to give an account, intended for the non-specialist, of the two proofs of the Göttsche conjecture given, respectively, first by Tzeng, and shortly after by Kool, Shende and Thomas.

Given a smooth, complex, projective surface $S$ and a positive integer $\delta$, the Göttsche conjecture predicts that if $L$ is a sufficiently positive line bundle on $S$, then the number of $\delta$-nodal curves in a general $\delta$-dimensional linear subsystem $\Lambda \subseteq |L|$ is given by a universal polynomial of degree $\delta$ in the four numbers $L^2$, $LK_S$, $K_S^2$, and $c_2(S)$.

Besides, the Göttsche–Yau–Zaslow formula expresses the generating function of the aforementioned numbers of $\delta$-nodal curves in terms of three explicit quasimodular forms and two unknown universal power series, whose coefficients can be determined by recursion.

## 1 – Introduction

The main theme of this note is the set of enumerative conjectures stated and discussed by Göttsche in [8]. We retrace and summarize the proofs of the Göttsche conjecture given first by Tzeng [33], and shortly after by Kool, Shende and Thomas [18]. We also discuss the proof of the Göttsche–Yau–Zaslow formula given by Tzeng. Accordingly, [8], [18] and [33] are the main sources for this note, and the proofs given here are directly inspired by those included in the original papers.

Let $S$ be a smooth complex projective surface, and let $\delta$ be a positive integer. A *$\delta$-nodal* curve on $S$ is a reduced (possibly reducible) curve $C \subseteq S$ having exactly $\delta$ nodes (i.e., ordinary double points), and no other singularities. Given a line bundle $L$ on $S$, we denote by $a_\delta(S, L)$

This text is from a series originated in the seminar *Degenerations and enumeration of curves on surfaces*, held at the University of Roma "Tor Vergata" in the years 2015–2017. See `https://www.math.univ-toulouse.fr/~tdedieu/#enumTV` for the full collection.

the number of $\delta$-nodal curves contained in a general linear subsystem $\Lambda \subseteq |L|$ of dimension $\delta$. When $L$ is sufficiently ample, $\delta$-nodal curves in $|L|$ indeed occur in codimension $\delta$, so that $a_\delta(S, L)$ is a finite number; moreover, it coincides with the number of $\delta$-nodal curves passing through $\dim |L| - \delta$ general points of $S$, see [16, p.234], reproduced in [XIII, Lemma **]. The study of $a_\delta(S, L)$ is a very classical theme, and the main topic of this note. We refer the reader to [17, 18, 33], and [XIII] in this volume, for a detailed bibliography and a treatise of the background underlying the results we are going to discuss.

The first result we would like to present is the so-called Göttsche conjecture [8, Conjecture 2.1]. In the light of [34, 16] (described in [XIII]) where formulas for the numbers $a_\delta(S, L)$ were given for all $\delta \leqslant 8$, it predicts that if the line bundle $L$ on $S$ is sufficiently ample with respect to $\delta$, then the number $a_\delta(S, L)$ is computed by a universal polynomial of degree $\delta$ in the numbers $L^2, LK_S, K_S^2, c_2(L)$ (universal means that the same polynomial works for all pairs $(S, L)$). Göttsche's conjecture has been proved in terms of the following positivity notion for a line bundle $L$ on $S$.

**(1.1) Definition.** *Let $L$ be a line bundle on a smooth projective surface $S$, and let $k \geqslant 0$ be an integer. We say that $L$ is* $k$-very ample *if for any $0$-dimensional subscheme $Z \subseteq S$ of length $k+1$, the restriction map $H^0(S, L) \to H^0(Z, L \otimes \mathcal{O}_Z)$ is surjective.*

The Göttsche conjecture was first proved by Tzeng [33]. Shortly after, Kool, Shende and Thomas proved the following version (see [18, Theorem 4.1]), with weaker positivity assumptions on $L$.

**(1.2) Theorem** (Göttsche conjecture)**.** *For all positive integers $\delta$, there exists a polynomial $T_\delta(x, y, z, t)$ of degree $\delta$ with the following property: for any smooth complex projective surface $S$ and for any $\delta$-very ample line bundle $L$ on $S$, a general $\delta$-dimensional linear subsystem $\Lambda \subseteq |L|$ contains exactly $T_\delta(L^2, LK_S, K_S^2, c_2(S))$ $\delta$-nodal curves, i.e.,*

$$a_\delta(S, L) = T_\delta(L^2, LK_S, K_S^2, c_2(S)). \tag{1.2.1}$$

It is worth noting that the coefficients of the polynomial $T_\delta(x, y, z, t)$ may be computed for all $\delta$. Indeed, the numbers $a_\delta(S, L)$ can be determined for all line bundles on $\mathbf{P}^2$ by [3, Theorem 1.1], and for all primitive line bundles on a $K3$ surface by [2, Theorem 1.1]; moreover, the 4-tuple of topological invariants $(L^2, LK_S, K_S^2, c_2(L))$ is $(n^2, -3n, 9, 3)$ if $(S, L) = \big(\mathbf{P}^2, \mathcal{O}_{\mathbf{P}^2}(n)\big)$, and $(L^2, 0, 0, 24)$ if $L$ is a primitive line bundle on a $K3$ surface $S$; it follows that the coefficients of each $T_\delta(x, y, z, t)$ can be computed by solving a system of linear equations obtained from (1.2.1) using different values of $n$ and $L^2$.

We will discuss the two proofs of Theorem (1.2) given by Tzeng [33], and Kool, Shende, and Thomas [18], respectively. There is also another proof due to Liu [27, 28], but we will not report on it. The Kool–Shende–Thomas proof is close in spirit to the work of Vainsencher and Kleiman–Piene on node polynomials, which is described in [XIII]. In particular it enables one to compute explicitly the node polynomials $T_\delta$ although, in its authors' own words, it is not terribly efficient. Tzeng's proof, on the other hand, involves an abstract version of the degeneration methods constituting the backbone of this volume, see, in particular, [I] and [VI].

We retrace the Kool–Shende–Thomas proof in Section 2. It relies on BPS calculus on Hilbert schemes of points on (arbitrary reduced) curves on a surface $S$, see [30, 32], the calculation of tautological integrals on Hilbert schemes of points on $S$ by means of the recursion introduced in [6], and the deformation theory of plane singularities of curves, see, e.g., [4] and [II]. The main idea underlying the proof is to describe the numbers $a_\delta(S, L)$ of $\delta$-nodal curves in a general linear subsystem $\Lambda \subseteq |L|$ in terms of the Euler characteristics of relative Hilbert schemes of points associated to the universal curve $\mathcal{C} \to \Lambda$. This is what BPS calculus is used for; after that, the proof consists in (non-trivial) intersection theoretic considerations. We refer to the introduction of Section 2 for a more detailed outline.

We survey Tzeng's proof in Section 3. Ultimately, it relies on the fact, proved in [20], that the algebraic cobordism group of smooth, complex, projective surfaces is generated by (the classes of) $\mathbf{P}^2$ and $\mathbf{P}^1 \times \mathbf{P}^1$. Tzeng's proof is formulated in the framework of [21, 22] for algebraic cobordism, which is based on the notion of (extended) double point relation: roughly speaking, one mods out the free abelian group generated by all smooth surfaces by relations corresponding to degenerations of one surface to the transverse union of two surfaces. Tzeng's main result is a suitable version of Jun Li's celebrated degeneration formula in this context, which involves Li and Wu's construction of a moduli space of stable ideal sheaves in [25]. The degeneration formula enables her to express the numbers $a_\delta$ for the pair $(S, L)$ (actually, the better suited variant $d_\delta(S, L)$ introduced by Göttsche) in terms of the numbers $a_\delta$ for the generators of the algebraic cobordism group, which then leads to the Göttsche conjecture.

The other results we will discuss concern the structure of the generating function of the universal polynomials $T_\delta(L^2, LK_S, K_S^2, c_2(S))$, that is the following formal power series in $x$,

$$T(S,L) := \sum_{\delta \geqslant 0} T_\delta(L^2, LK_S, K_S^2, c_2(S))\, x^\delta$$

(we set $T_0(x, y, z, t) = 1$). In [8] Göttsche proved that if the numbers $a_\delta(S, L)$ are indeed computed by universal polynomials, then the function $T(S, L)$ must satisfy rather strong conditions. In particular, it has the following multiplicative structure by [8, Proposition 2.3], which we may now state unconditionally since the Göttsche conjecture is proved.

**(1.3) Theorem.** *There exist universal power series $A_1, A_2, A_3, A_4 \in \mathbf{Q}[[x]]^\times$ such that*

$$T(S,L) = A_1^{L^2} A_2^{LK_S} A_3^{K_S^2} A_4^{c_2(S)}.$$

We emphasize that the power series $A_i$ are invertible in $\mathbf{Q}[[x]]$ (in fact, their constant coefficients all equal 1), so their exponentiations are defined unambiguously.

Göttsche's proof of this fact consists in considering a pair $(S, L)$ obtained as the disjoint union of pairs of the form $(\mathbf{P}^2, \mathcal{O}_{\mathbf{P}^2}(n))$ and $(\mathbf{P}^1 \times \mathbf{P}^1, \mathcal{O}_{\mathbf{P}^1\times\mathbf{P}^1}(n, n))$, and using the additivity properties of the 4-tuples $(L^2, LK_S, K_S^2, c_2(S))$. We will give Tzeng's proof in Section 3.4; it formulates the same idea in the somehow taylor-made theory of algebraic cobordism.

We further point out that, as was the case for the coefficients of the universal polynomials $T_\delta(x, y, z, t)$, the coefficients of the universal power series $A_1, A_2, A_3, A_4$ may be determined by the recursive formulas describing the numbers $a_\delta(S, L)$ for certain classes of pairs $(S, L)$ (see, e.g., the aforementioned [3, Theorem 1.1] for line bundles on $\mathbf{P}^2$ and [35, Theorem 6.7] for the case of Hirzebruch surfaces).

Another important result concerning the generating function $T(S, L)$ is the Göttsche–Yau–Zaslow formula, which was conjectured by Göttsche [8] based on the formulas in [34, 16], as a generalization of the Yau–Zaslow formula for rational curves on $K3$ surfaces [37], see [XI]. In particular, the Göttsche–Yau–Zaslow formula expresses the universal power series $T(S, L)$ in terms of five universal generating functions in one variable $q$: three explicit quasimodular forms, and two unknown universal power series, whose coefficients may nevertheless be determined by the recursion formula of [3].

We refer the reader to [14] for some preliminary notions on modular and quasimodular forms. Let $\tau \in \mathbb{C}$ be a complex number with positive imaginary part, and set $q := e^{2\pi i\tau}$. Consider the *discriminant form*

$$\Delta(\tau) := q \prod_{n>0} (1 - q^n)^{24}$$

and the *second Eisenstein series*

$$G_2(\tau) := -\frac{1}{24} + \sum_{n>0}\Big(\sum_{k|n} k\Big)q^n,$$

which are a modular form and a quasimodular form, respectively. Furthermore, given the differential operator $D := \frac{1}{2\pi i}\frac{d}{d\tau} = q\frac{d}{dq}$, we have that for any quasimodular form $f$, the derivative $Df$ is again a quasimodular form. Then the following holds, see [8, Conjecture 2.4] and [33, Theorem 1.2].

**(1.4) Theorem** (Göttsche-Yau–Zaslow formula)**.** *There exist universal power series $B_1(q)$ and $B_2(q)$ such that*

$$\sum_{\delta\geqslant 0} T_\delta(L^2, LK_S, K_S^2, c_2(S)) \cdot DG_2(\tau)^\delta = \frac{(DG_2(\tau)/q)^{\chi(L)}\, B_1(q)^{K_S^2} B_2(q)^{LK_S}}{(\Delta(\tau)D^2G_2(\tau)/q^2)^{\chi(\mathcal{O}_S)/2}}.$$

We give the proof in Section 3.4. An essential input in it is Bryan and Leung's main result in [2], see [XI].

We note that, according to [8, Remark 2.6] and [33, Corollary 4.4], the Göttsche–Yau–Zaslow formula can be reformulated as follows. For a smooth projective surface $S$ and for any $l, m, r \in \mathbf{Z}$, we define

$$N_r^S(l,m) := T_{l+\chi(\mathcal{O}_S)-1-r}(2l+m, m, K_S^2, c_2(S)).$$

Then Theorem (1.4) gives that

$$\sum_{l\in\mathbf{Z}} N_r^S(l,m)q^l = B_1(q)^{K_S^2} B_2(q)^m (DG_2(\tau))^r \frac{D^2G_2(\tau)}{(\Delta(\tau)D^2G_2(\tau)/q^2)^{\chi(\mathcal{O}_S)/2}},$$

and Theorem (1.2) ensures that if $L$ is a $\delta$-very ample line bundle on $S$ with $\delta = \chi(L) - 1 - r$, then the coefficient $N_r^S\left(\frac{L^2-LK_S}{2}, LK_S\right)$ counts the number of $\delta$-nodal curves in a general linear subsystem $W \subseteq |L|$ of codimension $r$. Analogously, setting $g = \frac{L^2+LK_S}{2} + 1 - \delta$, the coefficient

$$M_g^S\left(\tfrac{L^2-LK_S}{2}, LK_S\right) := N_{g-m-2+\chi(\mathcal{O}_S)}\left(\tfrac{L^2-LK_S}{2}, LK_S\right)$$

computes the number of nodal curves of geometric genus $g$ in a general linear subsystem of $|L|$ having codimension $r = g - LK_S - 2 + \chi(\mathcal{O}_S)$.

To conlude this introduction, let us mention some further developments, which we will not discuss any more in the text. Using the BPS calculs described in Section 2.2, Göttsche and Shende [9], see also [10], defined refined counting invariants, which interpolate between the usual counts over the complex numbers and Welschinger invariants; this, in turn, influenced tropical counts by Block-Goettsche [1].

Kleiman and Piene have further conjectured that Theorem (1.2) should hold more generally for a relative family of curves in a family of surfaces: this has been proved by Laarakker [19], by elaborating on the proof of Kool–Shende–Thomas. Rennemo [31] has furthermore proved a generalization of Theorem (1.2) valid for all isolated singularity types and in all dimensions, again following the way opened by Kool–Shende–Thomas.

**Acknowledgements.** The first author would like to thank Giulio Codogni, Lothar Göttsche, Manfred Lehn, Donatella Iacono, Margarida Melo, and Filippo Viviani, for helpful discussions and explanations during the preparation of a series of lectures on this topic and during the writing of this note.

The second author thanks Laurent Manivel and Adrien Sauvaget for their help with the Ellingsrud–Göttsche–Lehn recursion, and the first author for everything.

We both thank Martijn Kool and Richard Thomas for their comments on this text.

**(1.5) Notation.** We work throughout over the field $\mathbf{C}$ of complex numbers. By *curve* we mean a connected complete reduced algebraic curve, unless otherwise stated.

Given a projective variety $X$, we say that a property holds for a *general* point $x \in X$ if it holds on an open non-empty subset of $X$. In particular, if $L$ is a line bundle on a smooth surface $S$, by *general* linear subsystem $\Lambda \subseteq |L|$ of dimension $\delta$ we mean that $\Lambda$ is parameterized over an open subset of the Grassmannian of $\delta$-dimensional linear subspaces of $|L| \simeq \mathbf{P}^{h^0(L)-1}$.

For a positive integer $k$, we denote by $X^{[k]}$ the $k$-fold Hilbert scheme of points on a variety $X$, parameterizing 0-dimensional subschemes of $S$ of length $k$, and by $\mathrm{Hilb}^k(\mathcal{X}/B)$ the $k$-fold relative Hilbert scheme of points on the fibres of a family $\mathcal{X} \to B$ of varieties.

The *Euler number* (sometimes Euler characteristic) of a topological space $X$ is the number $e(X) = \sum_i (-1)^i \dim H^i_c(X, \mathbf{Z})$, where the subscript $c$ means cohomology with compact support. We refer to [XI, Section **] for some basic properties of the Euler number, including its additivity, and an expression when $X$ is a fibration.

# 2 – Kool–Shende–Thomas' proof of the Göttsche conjecture

In this section, we present the proof of Theorem (1.2) given by Kool, Shende, and Thomas [18]. So, we consider a line bundle $L$ on a smooth surface $S$, and we assume that $L$ is $\delta$-very ample for some integer $\delta \geqslant 1$.

In §2.1, we show that a general, $\delta$-dimensional, linear subsystem $\Lambda$ of $|L|$ contains a finite number of $\delta$-nodal curves, each occurring with multiplicity 1, and that all the other elements of $\Lambda$ are reduced curves having larger geometric genus. This is an application of the deformation theory of plane curve singularities, using in an essential way the $\delta$-very ampleness of $L$.

In §2.2, we first review some facts concerning the generating functions of the Euler characteristics $e(C^{[k]})$ of Hilbert schemes of points on curves $C$, for any reduced curve $C$ contained in $S$; this is the so-called BPS calculus on Hilbert schemes of points. We then use these facts to express the number $a_\delta(S, L)$ of $\delta$-nodal curves in a general $\Lambda \subseteq |L|$ as a linear combination of the Euler numbers $e(\mathrm{Hilb}^k(\mathcal{C}/\Lambda))$, $k = 0, \ldots, \delta$, of the relative Hilbert schemes of points on the fibres of the universal curve $\mathcal{C} \to \Lambda$. This is done by stratifying $\Lambda$ according to the topological types of its members, and using the additivity of the Euler number.

In §2.3, we conclude the proof by showing that each term appearing in this linear combination is indeed a polynomial in $L^2$, $LK_S$, $K_S^2$ and $c_2(L)$. To this aim, we retrace Kool–Shende–Thomas' computation of the Euler numbers $e(\mathrm{Hilb}^k(\mathcal{C}/\Lambda))$ in terms of certain tautological integrals, first on the Hilbert schemes $S^{[k]}$, and then on the products $S^k$.

The passing from integrals on $S^{[k]}$ to integrals on $S^k$ relies on the recursion of Ellingsrud, Göttsche and Lehn [6, Sections 3 and 4], which we survey in §2.4.

**(2.1)** We point out that the idea underlying this count of $\delta$-nodal curves generalizes the following argument for the case $\delta = 1$ (it is explained in more details in [XI, Section **]). Assume that $L$ is a very ample line bundle and $\Lambda \subseteq |L|$ is a general pencil, and denote by $p$ the genus of the smooth fibres of the universal curve $\mathcal{C} \to \Lambda \cong \mathbf{P}^1$.

The quantity $e(C) - (2 - 2p)$ equals 1 if $C$ is 1-nodal, and 0 if $C$ is smooth. Therefore, summing this quantity over $\Lambda \cong \mathbf{P}^1$, we obtain the identity

$$a_1(S, L) = e(\mathcal{C}) - (2 - 2p)\, e(\mathbf{P}^1).$$

To see this formula as a particular case of the Kool–Shende–Thomas' count, we rewrite it in terms of linear combinations of Euler characteristics of Hilbert schemes of points: one has

$$e(C) - (2 - 2p) = e(C^{[1]}) - (2 - 2p)\, e_0(C^{[0]})$$

and, correspondingly,

$$e(\mathcal{C}) - (2 - 2p)\, e(\mathbf{P}^1) = e(\mathrm{Hilb}^1(\mathcal{C}/\Lambda)) - (2 - 2p)\, e(\mathrm{Hilb}^0(\mathcal{C}/\Lambda)).$$

In §2.2, we give a similar expression of $a_\delta(S, L)$ as a sum over all members $C$ of $\Lambda \cong \mathbf{P}^\delta$ of some linear combination of the Euler characteristics $e(C^{[k]})$, $k = 1, \ldots, \delta$, which takes value 1 for $\delta$-nodal curves, and 0 for curves having higher geometric genus. So, in order to prove the Göttsche conjecture, it remains to express the corresponding linear combination of the Euler numbers $e(\mathrm{Hilb}^k(\mathcal{C}/\Lambda))$, $k = 0, \ldots, \delta$, in terms of the topological invariants $L^2$, $LK_S$, $K_S^2$ and $c_2(S)$; this is the purpose of §2.3.

For $\delta = 1$, this goes simply as follows: $e(\mathcal{C}) = c_2(S) + L^2$ as $\mathcal{C}$ is the blow-up of $S$ at $L^2$ points, $2p - 2 = L(L + K_S)$ by the adjunction formula, and $e(\mathbf{P}^1) = 2$, so that we obtain

$$a_1(S, L) = c_2(S) + 3L^2 + 2LK_S,$$

which is indeed a polynomial of degree 1 in the topological numbers $L^2, LK_S, K_S^2, c_2(S)$.

## 2.1 – Sufficiently ample line bundles

Let $L$ be a line bundle on a smooth projective surface $S$. We follow [18, Section 2] to show how the $\delta$-very ampleness condition on $L$ ensures that the singularities of curves in a general dimension $\delta$ linear subsystem $\Lambda \subseteq |L|$ behave sufficiently well for the count of $\delta$-nodal curves in $\Lambda$ by means of the Euler numbers $e(\mathrm{Hilb}^k(\mathcal{C}/\Lambda))$ in subsection 2.2 to make sense, see Theorem (2.10). Roughly speaking, the singularities of the members of $\Lambda$ behave as expected up to codimension $\delta$ in $|L|$. The precise statement is the following.

**(2.2) Theorem.** *Assume $L$ is $\delta$-very ample for some positive integer $\delta$. Let $p$ denote the common arithmetic genus of all members of $|L|$, and consider a general linear subsystem $\Lambda \subseteq |L|$ of dimension $\delta$. Then all members of $\Lambda$ are reduced curves, and have geometric genus $g \geqslant p - \delta$. Moreover, there are finitely many members of genus $g = p - \delta$; they are all $\delta$-nodal curves, and appear with multiplicity 1.*

Let us clarify the last sentence. Recall that the locus $V^\delta(\Lambda)$ parametrizing $\delta$-nodal members of $\Lambda$ is well-defined as a subscheme of $\Lambda$, see [5] for a discussion and the relevant references; the locus $V_g(\Lambda)$ parametrizing genus $g$ members of $\Lambda$, however, is only defined set-theoretically (it is a locally closed subset of $\Lambda$). Now, the last sentence means that (i) the set $V_{p-\delta}(\Lambda)$ has dimension 0 and coincides with the support of $V^\delta(\Lambda)$, and (ii) the scheme $V^\delta(\Lambda)$ is reduced.

**(2.3) Remark.** Note that the theorem implies that all members of $\Lambda$ have at most $\delta$ singular points, and all members with $\delta$ singular points are $\delta$-nodal curves. Indeed, each singular point makes the geometric genus drop by at least one.

The rest of Section 2.1 is dedicated to the proof of Theorem (2.2). We will use extensively the deformation theory of planar curve singularities, for which we refer to [II] in this volume, and to [4, 11] for extended treatises; the reader may also consult [5]. The punchline is that the $\delta$-very ampleness property ensures that the maps from the relevant equigeneric and equisingular subloci of $|L|$ to the corresponding semi-universal deformation spaces of singularities are smooth, because their differentials are restriction maps, and the $\delta$-very ampleness grants their surjectivity.

**(2.4)** Let us first show that the subscheme $V^\delta(|L|) \subseteq |L|$ of $\delta$-nodal curves in $|L|$ is smooth of codimension $\delta$. By Bertini's theorem this implies that $V^\delta(\Lambda) = V^\delta(|L|) \cap \Lambda$ is smooth of dimension 0, i.e., $\delta$-nodal curves in $\Lambda$ are finitely many and appear with multiplicity 1.

Let $[C] \in |L|$ be a $\delta$-nodal curve, and $J$ be its Jacobian ideal. The vector space $B := H^0\left(C, \mathcal{O}_C/J\right)$ is the product of the semi-universal deformation spaces of the various singularities of $C$. We thus have a map of pointed germs $\varphi : (|L|, [C]) \to (B, 0)$, and its differential is the restriction map

$$\rho_J : T_{|L|,[C]} \cong H^0(S, L)/\langle s_C \rangle \longrightarrow H^0\left(Z, L \otimes \mathcal{O}_C/J\right),$$

where $Z \subseteq S$ is the singular subscheme of $C$, defined in $C$ by the ideal $J$, and $s_C \in H^0(S, L)$ is a section with zero divisor equal to $C$. Since $C$ is $\delta$-nodal, its singular subscheme has length $\delta$, hence the restriction map $H^0(S, L) \to H^0(Z, L \otimes \mathcal{O}_Z)$ is surjective by the positivity property of $L$. This implies that the map $\rho_J$ is surjective, hence $\varphi$ is smooth.

The semi-universal deformation space of a node is a line, in which the equisingular locus is only the origin. Thus the locus of $\delta$-nodal curves in $B$ is only the origin, which is smooth and has codimension $\delta$ in $B$. By smoothness of $\varphi$, the same holds for $V^\delta(|L|)$ around $[C]$, which proves our claim.

**(2.5)** Consider a non-negative integer $g' < p - \delta$, and let $V_{g'}(|L|) \subseteq |L|$ be the locally closed subset of reduced curves having geometric genus $g'$. We now show that $V_{g'}(|L|)$ has codimension (strictly) greater than $\delta$ in $|L|$, which implies that $\Lambda$ does not contain any curve of genus $g'$.

Let $[C] \in V_{g'}(|L|)$. We consider its conductor ideal $A$, and let $Y$ be the subscheme of $C$ defined by $A$. By [II, Proposition 3.1] the tangent cone of $V_{g'}(|L|)$ at $[C]$ is contained in the kernel of the restriction map

$$\rho_Y : H^0(S, L)/\langle s_C \rangle \longrightarrow H^0\left(Y, L \otimes \mathcal{O}_C/A\right),$$

where $s_C$ is again a section with zero divisor equal to $C$.

Since the co-genus $p - g'$ of $C$ in $|L|$ is greater than $\delta$, we may find a subscheme $Z \subseteq Y$ of length $\delta + 1$. Then the positivity property of $L$ implies that the restriction map

$$\rho_Z : H^0(S, L)/\langle s_C \rangle \longrightarrow H^0\left(Z, L \otimes \mathcal{O}_Z\right)$$

is surjective, hence its kernel has codimension $\delta + 1$. Since $\rho_Z = \pi \circ \rho_Y$, where $\pi$ is the restriction map $H^0(L \otimes \mathcal{O}_Y) \to H^0(L \otimes \mathcal{O}_Z)$, this implies that the kernel of $\rho_Y$ has codimension at least $\delta + 1$, hence $V_{g'}(|L|)$ has codimension at least $\delta + 1$ in $|L|$ as required.

**(2.6)** Let $\tau$ be a collection of topological types of planar singularities, which does not consist of $\delta$ nodes, and such that a member of $|L|$ with singularities of type $\tau$ has genus $p - \delta$. We now show that the subscheme $V^\tau(|L|) \subseteq |L|$ parametrizing curves with topological type $\tau$ has codimension (strictly) greater than $\delta$ in $|L|$.

Let $[C] \in V^\tau(|L|)$, and $I$ be its equisingular ideal. By [II, Proposition 3.1], the tangent space to $V^\tau(|L|)$ at $[C]$ is the kernel of the restriction map

$$\rho_Y : H^0(S, L)/\langle s_C \rangle \longrightarrow H^0\left(Y, L \otimes \mathcal{O}_C/I\right),$$

where $Y$ is the subscheme of $C$ defined by $I$. The assumptions on $\tau$ imply that we may find a subscheme $Z \subseteq Y$ of length $\delta+1$, see [II, Lemma 2.4]. Then, as in (2.5), the $\delta$-very ampleness of $L$ implies that the kernel of $\rho_Y$ has codimension at least $\delta+1$, hence $V^\tau(|L|)$ has codimension at least $\delta+1$ in $|L|$ as claimed.

**(2.7)** Finally, we proceed to show that the locus $V^{\mathrm{nr}}(|L|) \subseteq |L|$ parametrizing non-reduced curves has codimension (strictly) greater than $\delta$.

We claim that if $C$ is a non-reduced member of $|L|$, then the tangent cone to $V^{\mathrm{nr}}(|L|)$ at $[C]$ is contained in the kernel of the restriction map

$$\rho : H^0(S,L)/\langle s_C \rangle \longrightarrow H^0\left(C_{\mathrm{red}}, L \otimes \mathcal{O}_{C_{\mathrm{red}}}\right).$$

When we know this, we conclude that $V^{\mathrm{nr}}(|L|)$ has codimension greater than $\delta$ always along the same lines: we pick a 0-dimensional subscheme $Z \subseteq C_{\mathrm{red}}$ of length $\delta+1$, and the kernel of $\rho$ is contained in the kernel of the restriction to $Z$, which has codimension $\delta+1$ since $L$ is $\delta$-very ample.

To prove our claim, we argue as in [36, §2.1] and consider the incidence variety

$$\mathcal{D} := \left\{ \big([C],p\big) \in |L| \times S : C \text{ is singular at } p \right\}.$$

It is smooth: indeed since $L$ is very ample and $S$ is smooth, $\mathcal{D}$ is a projective bundle over $S$. On the other hand, if $[C] \in V^{\mathrm{nr}}(|L|)$ then $([C],p) \in \mathcal{D}$ for all $p \in C_{\mathrm{red}}$. It follows that the tangent cone to $V^{\mathrm{nr}}(|L|)$ at $[C]$ is contained in the image of the tangent space to $\mathcal{D}$ at $([C],p)$ by the projection $T_{|L|\times S} \to T_{|L|}$ for all $p \in C_{\mathrm{red}}$. We shall show that the latter image is contained in the kernel of the restriction

$$H^0(S,L)/\langle s_C \rangle \longrightarrow H^0\big(p, L|_p\big),$$

and this will prove our claim.

The equations of $\mathcal{D}$ in $|L| \times S$ are as follows. Let $([C],p) \in |L| \times S$; we choose a local system of coordinates $(x,y)$ on $S$ around $p$, and see $|L|$ around $[C]$ as a finite dimensional affine space of polynomials $f \in \mathbf{C}[x,y]$. In this representation, $\mathcal{D}$ is defined by the equations

$$f(x,y) = df(x,y) = 0.$$

By differentiation, the latter equations give the equations of the tangent space of $\mathcal{D}$ in $T_{|L|} \times T_S$. In particular, if $(h,v) \in T_{|L|,f} \times T_{S,(x,y)}$ sits in $T_{\mathcal{D}}$, then

$$df(x,y)(v) + h(x,y) = 0,$$

which is equivalent to $h(x,y) = 0$ since $df(x,y) = 0$. In other words, if $([C],p) \in \mathcal{D}$, then the projection of $T_{\mathcal{D},([C],p)}$ to $T_{|L|,[C]}$ is contained in the kernel of the restriction map

$$H^0(S,L)/\langle s_C \rangle \longrightarrow H^0\big(p, L|_p\big),$$

as we wanted to prove. We have thus proved our claim, hence that $V^{\mathrm{nr}}(|L|)$ has codimension greater than $\delta$.

The results proved in (2.4), (2.5), (2.6), and (2.7) together amount to Theorem (2.2). The proof of the latter is therefore complete. □

## 2.2 – Euler characteristics of Hilbert schemes of points over curves

We consider a smooth projective surface $S$ endowed with a $\delta$-very ample line bundle $L$ on $S$, as in Theorem (2.2). We follow [18, Section 3] in order to express the number $a_\delta(S,L)$ of $\delta$-nodal curves in a general $\delta$-dimensional linear subsystem $\Lambda \subseteq |L|$ in terms of Euler numbers $e(\mathrm{Hilb}^k(\mathcal{C}/\Lambda))$ of relative Hilbert schemes, where $\mathcal{C} \to \Lambda$ is the universal curve.

To this aim, we recall a couple of preliminary results describing the generating series of the Euler characteristics $e(C^{[k]})$ of $k$-fold Hilbert schemes of points on curves $C \subseteq S$, see [18, Propositions 3.1 and 3.2]. Both results were first proved in [30] in the setting of stable pairs, whereas [32, Section 2] includes proofs using only Hilbert schemes.

**(2.8) Proposition.** *Let $S$ be a smooth surface and let $C \subseteq S$ be a curve having arithmetic genus $p$ and geometric genus $g$. There exists a unique sequence of integers $(n_{C,g}, n_{C,g+1}, \dots, n_{C,p})$, depending only on the topological type of $C$, such that the generating series of the Euler characteristics $e(C^{[k]})$ satisfies*

$$\sum_{k=0}^{\infty} e(C^{[k]}) q^k = \sum_{i=g}^{p} n_{C,i} q^{p-i} (1-q)^{2i-2}.$$

Writing the generating series $\sum e(C^{[k]})q^k$ as a sum $\sum n_{C,i} q^{p-i}(1-q)^{2i-2}$ with $i$ ranging between $-\infty$ and $p$ is granted by a formal property of power series, see [32, Proposition 1]. However, the geometry of $C$ is involved in the proof that $n_{C,i} = 0$ for all $i < g$; this is deduced from [32, Proposition 8], see [32, Corollary 11].

When the curve $C$ is $\delta$-nodal, the generating series of Euler characteristics $e(C^{[k]})$ is governed by the following, see [32, Corollary 12].

**(2.9) Proposition.** *Let $C$ be a $\delta$-nodal curve with arithmetic genus $p$. Then the coefficients $n_{C,i}$ appearing in Proposition (2.8) are*

$$n_{C,i} = \binom{\delta}{p-i}$$

*for all $i = p-\delta, \dots, p$.*

When $C$ is a smooth curve this was proved by Macdonald [29], in order to give a new proof of the de Jonquières formulas discussed in [XII]. In this case one finds

$$\sum_{k=0}^{\infty} e\big(C^{[k]}\big) q^k = (1-q)^{2p-2} = \sum_{k=0}^{2p-2} (-1)^k \binom{2p-2}{k} q^k.$$

While $e(C) = 2-2p$, and $e(C^{[k]}) = 0$ if $k > 2p-2$ because then $C^{[k]} \to J^k(C)$ gives a structure of projective bundle over a torus, the values of $e(C^{[k]})$ for $1 < k \leqslant 2p-2$ are not obvious to us.

In the general case of Proposition (2.9), Kool, Shende, and Thomas interpret the number $\binom{\delta}{p-i}$ as the number of partial normalizations of $C$ at $p-i$ of its $\delta$ nodes. This is the reason why Proposition (2.9) is considered a statement in BPS calculus, see [XI].

We may now state the main result of this section, see [18, Theorem 3.4]. It follows from Propositions (2.8), (2.9) and Theorem (2.2)

**(2.10) Theorem.** *Let $S$ be a smooth projective surface equipped with a line bundle $L$, and denote by $p$ be the arithmetic genus of all members of $|L|$. Let $\delta$ be a positive integer such that $L$ is $\delta$-very ample, and consider a general linear subsystem $\Lambda \subseteq |L|$ of dimension $\delta$. Let $\mathcal{C} \to \Lambda$ be* *the universal curve over* $\Lambda$*. There are integers* $n_{\Lambda,p-\delta},\dots,n_{\Lambda,\delta}$ *such that the following identity of formal power series holds:*

$$\sum_{k=0}^{\infty} e(\mathrm{Hilb}^k(\mathcal{C}/\Lambda))q^k = \sum_{i=p-\delta}^{p} n_{\Lambda,i}q^{p-i}(1-q)^{2i-2}. \tag{2.10.1}$$

*Moreover, the number* $a_\delta(S,L)$ *of* $\delta$*-nodal curves in* $\Lambda$ *equals* $n_{\Lambda,p-\delta}$*.*

*Proof.* Let $T$ be the (finite) set indexing the topological types of the members of $\Lambda$. For all $\tau\in T$, let $V^\tau_{\mathrm{red}} = V^\tau(\Lambda)_{\mathrm{red}}$ be the locus of curves having topological type $\tau$, and let $g_\tau$ denote the common geometric genus of these curves. Note that $g_\tau \geqslant p-\delta$ by Theorem (2.2). Let $\mathcal{C}_\tau \to V^\tau_{\mathrm{red}}$ be the universal curve.

For all $k\geqslant 0$, the relative Hilbert scheme $\mathrm{Hilb}^k(\mathcal{C}_\tau/V^\tau_{\mathrm{red}})$ is a topological fibre bundle, so $e(\mathrm{Hilb}^k(\mathcal{C}_\tau/V^\tau_{\mathrm{red}})) = e(C_\tau^{[k]})e(V^\tau_{\mathrm{red}})$ for any membre $C_\tau$ of $\Lambda$ of topological type $\tau$. For all $i=g_\tau,\dots,p$, we denote by $n_{\tau,i}$ the number $n_{C_\tau,i}$ as in Proposition (2.8) (it depends only on $\tau$); we set $n_{\tau,i}=0$ if $i<g_\tau$. Then, by the additivity of the Euler number, and using Proposition (2.8), we have

$$\begin{aligned}\sum_{k=0}^{\infty} e(\mathrm{Hilb}^k(\mathcal{C}/\Lambda))q^k &= \sum_{k=0}^{\infty}\Big(\sum_{\tau\in T} e(C_\tau^{[k]})e(V^\tau_{\mathrm{red}})\Big)q^k = \sum_{\tau\in T}\Big(\sum_{k=0}^{\infty} e(C_\tau^{[k]})q^k\Big)e(V^\tau_{\mathrm{red}}) = \\ &= \sum_{\tau\in T}\Big(\sum_{i=g_\tau}^{p} n_{\tau,i}q^{p-i}(1-q)^{2i-2}\Big)e(V^\tau_{\mathrm{red}}) = \sum_{i=p-\delta}^{p} n_{\Lambda,i}q^{p-i}(1-q)^{2i-2},\end{aligned}$$

with $n_{\Lambda,i} := \sum_{\tau\in T} n_{\tau,i}e(V^\tau_{\mathrm{red}})$: this proves the identity (2.10.1).

Theorem (2.2) ensures that $\delta$-nodal curves are the only curves in $\Lambda$ having geometric genus $p-\delta$, they are finitely many, and they appear with multiplicity 1. Therefore $n_{\Lambda,p-\delta} = a_\delta(S,L)\cdot n_{\delta A_1,p-\delta}$, where $\delta A_1$ denotes the topological type of $\delta$-nodal curves members of $|L|$. Since $n_{\delta A_1,p-\delta}=1$ by Proposition (2.9), one has $n_{\Lambda,p-\delta} = a_\delta(S,L)$ as we wanted. □

**(2.11) Corollary.** *The number* $a_\delta(S,L)$ *is a linear combination of the numbers* $e(\mathrm{Hilb}^k(\mathcal{C}/\Lambda))$*,* $k=0,\dots,\delta$*, in which the coefficient of* $e(\mathrm{Hilb}^k(\mathcal{C}/\Lambda))$ *is a polynomial in* $p$ *of degree at most* $\delta-k$ *for all* $k$*, and the coefficient of* $e(\mathrm{Hilb}^\delta(\mathcal{C}/\Lambda))$ *is* 1*.*

*Proof.* Let $e_{\Lambda,k} = e(\mathrm{Hilb}^k(\mathcal{C}/\Lambda))$ for all $k=0,\dots,\delta$. Expanding $(1-q)^{2i-2}$ for all $i$ in (2.10.1) gives

$$\sum_{k=0}^{\infty} e_{\Lambda,k}q^k = \sum_{i=p-\delta}^{p} n_{\Lambda,i}\sum_{m=0}^{2i-2}\tbinom{2i-2}{m}(-1)^m q^{p-i+m}.$$

Then, by comparing the coefficients of $q^k$ for $k=0,\dots,\delta$, one sees with standard algebraic manipulations that

$$e_{\Lambda,k} = \sum_{m=0}^{k}(-1)^m\tbinom{2(p-k+m)-2}{m}n_{\Lambda,p-k+m}. \tag{2.11.1}$$

Let us now show by induction on $k=0,\dots,\delta$ that

$$n_{\Lambda,p-k} = P_{k,0}(p)\,e_{\Lambda,0} + \cdots + P_{k,k-1}(p)\,e_{\Lambda,k-1} + e_{\Lambda,k},$$

where each $P_{k,s}$ is a polynomial of degree at most $k-s$, independent of $\Lambda$. The case $k=\delta$ gives the wanted result, since $a_\delta(S,L) = n_{\Lambda,p-\delta}$ by Theorem (2.10).

For $k=0$ one has $n_{\Lambda,p}=e_{\Lambda,0}$ by (2.11.1), as required. Suppose now that $k$ is positive. For all $m=1,\dots,k$, we set $Q_{k,m}(p)=(-1)^m\binom{2(p-k+m)-2}{m}$, a polynomial of degree $m$ in $p$. By (2.11.1) and our recursion hypothesis, we have

$$e_{\Lambda,k}=n_{\Lambda,p-k}+\sum_{m=1}^{k}Q_{k,m}(p)\sum_{s=0}^{k-m}P_{k-m,s}(p)e_{\Lambda,s}=n_{\Lambda,p-k}+\sum_{s=0}^{k-1}\sum_{m=1}^{k-s}Q_{k,m}(p)\,P_{k-m,s}(p)e_{\Lambda,s}$$

(setting $P_{m,m}=1$). This gives an explicit expression for $P_{k,s}$, which ends the proof since

$$\deg(Q_{k,m}P_{k-m,s})\leqslant m+(k-m-s)=k-s.$$

□

## 2.3 – Proof of the Göttsche conjecture

In this section we prove Theorem (1.2), following [18, Section 4]. With Corollary (2.11) at hand, the remaining task is to compute the Euler characteristics $e(\mathrm{Hilb}^k(\mathcal{C}/\Lambda))$, expressing them as polynomials of degree $k$ in the numbers $L^2, LK_S, K_S^2, c_2(L)$; this is done with the recursion of [6].

*Proof of Theorem (1.2).* Consider a smooth projective surface $S$ and a $\delta$-very ample line bundle $L$ on $S$. Let $p$ be the arithmetic genus of the curves in $|L|$, and let $\Lambda\subseteq|L|$ be a general linear subsystem of dimension $\delta$ with universal curve $\mathcal{C}\to\Lambda$. Corollary (2.11) ensures that the number $a_\delta(S,L)$ of $\delta$-nodal curves in $\Lambda$ is a linear combination of the Euler numbers $e(\mathrm{Hilb}^k(\mathcal{C}/\Lambda))$, $k=0,\dots,\delta$, where the coefficient of each $e(\mathrm{Hilb}^k(\mathcal{C}/\Lambda))$ is a polynomial of degree at most $\delta-k$ in $p$, and the coefficient of $e(\mathrm{Hilb}^\delta(\mathcal{C}/\Lambda))$ equals 1. Moreover, these polynomials are independent of the pair $(S,L)$ and, since $p=\frac{1}{2}L^2+\frac{1}{2}K_SL+1$ by the adjunction formula, the coefficient of $e(\mathrm{Hilb}^k(\mathcal{C}/\Lambda))$ is a polynomial of degree at most $\delta-k$ in $L^2$ and $K_SL$. The Göttsche conjecture therefore follows from Proposition (2.12) below. □

**(2.12) Proposition.** *The Euler characteristics $e(\mathrm{Hilb}^k(\mathcal{C}/\Lambda))$ is a polynomial of degree exactly $k$ in the numbers $L^2, LK_S, K_S^2, c_2(L)$.*

*Proof.* The case $k=0$ follows from the following equalities:

$$e(\mathrm{Hilb}^0(\mathcal{C}/\Lambda))=e(\Lambda)=e(\mathbf{P}^\delta)=\delta+1.$$

Let us now assume that $1\leqslant k\leqslant\delta$. We divide the proof in three steps.

*Step 1.* We consider the Hilbert scheme $S^{[k]}$ and the universal subscheme $Z_k\subseteq S\times S^{[k]}$, endowed with the projections $\pi_1\colon Z_k\to S$ and $\pi_2\colon Z_k\to S^{[k]}$. We define a vector bundle $L^{[k]}$ of rank $k$ on $S^{[k]}$,

$$L^{[k]}:=\pi_{2*}\pi_1^*L.$$

Letting $\mathbf{P}:=\mathbf{P}\left(H^0(S,L)\right)=|L|$, we consider the line bundle $L\boxtimes\mathcal{O}_{\mathbf{P}}(1):=\pi_S^*L\otimes\pi_{\mathbf{P}}^*\mathcal{O}_{\mathbf{P}}(1)$ on $S\times\mathbf{P}$, where $\pi_S$ and $\pi_{\mathbf{P}}$ are the two projections. In view of the standard identification $H^0(S,L)\otimes H^0(S,L)^*\cong\mathrm{Hom}\left(H^0(S,L),H^0(S,L)\right)$, we may consider the canonical section $\mathrm{id}_{H^0(S,L)}\in H^0(S,L)\otimes H^0(S,L)^*$ of the line bundle $L\boxtimes\mathcal{O}_{\mathbf{P}}(1)$. By taking the pullback of $\mathrm{id}_{H^0(S,L)}$ to $S\times S^{[k]}\times\mathbf{P}$, and pushing down to $S^{[k]}\times\mathbf{P}$ its restriction to $Z_k\times\mathbf{P}\subseteq S\times S^{[k]}\times\mathbf{P}$, we obtain a tautological section of the vector bundle $M^{[k]}$ on $S^{[k]}\times\mathbf{P}$ defined as

$$M^{[k]}=L^{[k]}\boxtimes\mathcal{O}_{\mathbf{P}}(1):=\pi_{S^{[k]}}^*L^{[k]}\otimes\pi_{\mathbf{P}}^*\mathcal{O}_{\mathbf{P}}(1).$$

By abuse of notation, we denote by $\mathcal{C} \to \mathbf{P}$ the universal curve over the complete linear system $|L|$. The relative Hilbert scheme $\mathrm{Hilb}^k(\mathcal{C}/\mathbf{P})$ is a smooth subvariety of $S^{[k]} \times \mathbf{P}$, of dimension $k + \dim(\mathbf{P})$ hence codimension $k$ in $S^{[k]} \times \mathbf{P}$. Since it coincides moreover with the vanishing locus of the tautological section of $M^{[k]}$ corresponding to $\mathrm{id}_{H^0(S,L)}$, its class in the Chow ring of $S^{[k]} \times \mathbf{P}$ is Poincaré dual to the Chern class $c_k(M^{[k]})$, that is $\left[\mathrm{Hilb}^k(\mathcal{C}/\mathbf{P})\right] = c_k(M^{[k]}) \cap \left[S^{[k]} \times \mathbf{P}\right]$ (see, e.g., [7, Example 3.2.16]). We note furthermore that the fibre of the projection $\mathrm{Hilb}^k(\mathcal{C}/\mathbf{P}) \to S^{[k]}$ over any 0-dimensional scheme $Z \in S^{[k]}$ is given by pairs $(Z, [s]) \in S^{[k]} \times \mathbf{P}$ such that $s$ vanishes at $Z$, i.e., $s \in \ker\big(H^0(S,L) \xrightarrow{\gamma_Z} H^0\left(Z, L \otimes \mathcal{O}_Z\right)\big)$. As $L$ is $(k-1)$-very ample for all $k = 1, \dots, \delta+1$, the map $\gamma_Z$ is surjective, so that the tautologial section of $M^{[k]}$ vanishes transversally with respect to the fibration $\mathrm{Hilb}^k(\mathcal{C}/\mathbf{P}) \to S^{[k]}$, whose fibres therefore have constant dimension. Moreover, since the linear subsystem $\Lambda \subseteq \mathbf{P}$ is general, the restriction $\mathrm{Hilb}^k(\mathcal{C}/\Lambda) \subseteq S^{[k]} \times \Lambda$ is smooth by Bertini's Theorem. Therefore, since the class $\left[\mathrm{Hilb}^k(\mathcal{C}/\Lambda)\right]$ is Poincaré dual in the Chow ring of $S^{[k]} \times \Lambda$ to the $k$-th Chern class of the restriction of $M^{[k]}$, the Euler characteristic of $\mathrm{Hilb}^k(\mathcal{C}/\Lambda)$ satisfies

$$\begin{aligned} e(\mathrm{Hilb}^k(\mathcal{C}/\Lambda)) &= \int_{\mathrm{Hilb}^k(\mathcal{C}/\Lambda)} c_k(T_{\mathrm{Hilb}^k(\mathcal{C}/\Lambda)}) = \int_{\mathrm{Hilb}^k(\mathcal{C}/\Lambda)} c_\bullet(T_{\mathrm{Hilb}^k(\mathcal{C}/\Lambda)}) \\ &= \int_{\mathrm{Hilb}^k(\mathcal{C}/\Lambda)} \frac{c_\bullet(T_{S^{[k]}\times\Lambda})}{c_\bullet(M^{[k]})} = \int_{S^{[k]}\times\Lambda} c_k(M^{[k]}) \frac{c_\bullet(T_{S^{[k]}\times\Lambda})}{c_\bullet(M^{[k]})}. \end{aligned}$$

*Step 2.* Let us now spell out the various terms in the last integral above. Let $\omega = c_1(O_\Lambda(1))$ be the hyperplane class on $\Lambda \simeq \mathbf{P}^\delta$. One has $c_\bullet(T_\Lambda) = \sum_{k=0}^{\delta} \binom{\delta+1}{k} \omega^k = (1+\omega)^{\delta+1}$ (see, e.g., [13, p.414]), hence $c_\bullet(T_{S^{[k]}\times\Lambda}) = c_\bullet(T_{S^{[k]}}) c_\bullet(T_\Lambda) = c_\bullet(T_{S^{[k]}})(1+\omega)^{\delta+1}$.

On the other hand, since $M^{[k]}$ is the tensor product of the rank $k$ vector bundle $\pi^*_{S^{[k]}} L^{[k]}$ with the line bundle $\pi^*_{\mathbf{P}} \mathcal{O}_{\mathbf{P}}(1)$, its restriction to $\Lambda$ satisfies $c_k(M^{[k]}) = \sum_{i=0}^k \omega^i c_{k-i}(L^{[k]})$ (see, e.g., [7, p. 55]). Finally, [7, Example 3.2.2] yields

$$\begin{aligned} c_\bullet(M^{[k]}) &= \sum_{p=0}^{k} c_p(M^{[k]}) = \sum_{p=0}^{k} \left( \sum_{i=0}^{p} \tbinom{k-i}{p-i} c_i(L^{[k]}) \omega^{p-i} \right) = \\ &= \sum_{i=0}^{k} c_i(L^{[k]}) \sum_{p=i}^{k} \tbinom{k-i}{p-i} \omega^{p-i} = \sum_{j=0}^{k} c_{k-j}(L^{[k]})(1+\omega)^j. \end{aligned}$$

Thus,

$$e(\mathrm{Hilb}^k(\mathcal{C}/\Lambda)) = \int_{S^{[k]}\times\Lambda} \frac{c_\bullet(T_{S^{[k]}})(1+\omega)^{\delta+1} \sum_{i=0}^k \omega^i c_{k-i}(L^{[k]})}{\sum_{j=0}^k c_{k-j}(L^{[k]})(1+\omega)^j} \tag{2.12.1}$$

$$= \int_{S^{[k]}} \mathrm{Coeff}_{\omega^\delta} \Big( \frac{c_\bullet(T_{S^{[k]}})(1+\omega)^{\delta+1} \sum_{i=0}^k \omega^i c_{k-i}(L^{[k]})}{\sum_{j=0}^k c_{k-j}(L^{[k]})(1+\omega)^j} \Big). \tag{2.12.2}$$

*Step 3.* We now end the proof with the recursion procedure developed in [6], see Proposition (2.14). The latter, applied $k$ times, turns the above integral over $S^{[k]}$ into the integral over $S^k$ of a polynomial in the Chern classes of $\mathrm{pr}_i^* L$, $\mathrm{pr}_i^* T_S$, and $\mathrm{pr}_{ij}^* T_\Delta$ $(i \neq j)$, where $\Delta \subseteq S \times S$ is the diagonal, $\mathrm{pr}_i : S^k \to S$ is the $i$-th projection, and $\mathrm{pr}_{ij} = (\mathrm{pr}_i, \mathrm{pr}_j)$. The latter integral results (with some non-trivial work, see, for instance, [6, Proposition 0.5 and its proof]) in a polynomial in $L^2$, $LK_S$, $K_S^2$, and $c_2(S)$, of degree at most $k$.

It remains to prove that this degree is exactly $k$. To do so, we concentrate on the $c_2(S)^k$ term. One looks the explicit formulas for the recursion straight in the eyes, viz. (2.14.3), (2.14.4), (2.14.5), (2.14.6), (2.14.7), and finds that the only contribution to the $c_2(S)^k$ term comes from the $c_0(L^{[k]})$ term. One may thus, for our present purpose, replace the fraction in (2.12.2) by

$$\frac{c_\bullet(T_{S^{[k]}})(1+\omega)^{\delta+1}\omega^k}{(1+\omega)^k} = c_\bullet(T_{S^{[k]}})\omega^k(1+\omega)^{\delta-k+1},$$

the coefficient of which in $\omega^\delta$ is $(\delta-k+1)\,c_\bullet(T_{S^{[k]}})$. The corresponding summand of the integral (2.12.2) is $(\delta-k+1)\,c_{2k}(T_{S^{[k]}})$, which contributes by $(\delta-k+1)/k!$ to the $c_2(S)^k$ term, the factorial coming from the $k$ successive applications of (2.14.2). □

**(2.13) Example.** Let us spell out this computation for the number of 1-nodal curves in a pencil, and see how it relates with the standard way to do it, which we recalled in the introduction to the present Section 2. In this case we have $\Lambda \simeq \mathbf{P}^1$. Equation (2.11.1) gives

$$e(\mathrm{Hilb}^1(\mathcal{C}/\Lambda)) = n_{\Lambda,p-1} - (2p-2)n_{\Lambda,p} \quad \text{and} \quad e(\mathrm{Hilb}^0(\mathcal{C}/\Lambda)) = n_{\Lambda,p}.$$

We have $e(\mathrm{Hilb}^0(\mathcal{C}/\Lambda)) = e(\Lambda) = 2$, and we concentrate on the computation of $e(\mathrm{Hilb}^1(\mathcal{C}/\Lambda))$ as the integral (2.12.2). The numerator of the integrand of (2.12.1) is

$$\begin{aligned}(1-K+c_2)(1+\omega)^2(L+\omega) &= (1-K+c_2)(1+2\omega)(L+\omega)\\ &= \big[L-KL\big] + \big[1-K+2L+c_2-2KL\big]\omega,\end{aligned}$$

and on the other hand

$$\begin{aligned}\frac{1}{1+L+\omega} &= 1-(L+\omega)+(L+\omega)^2-(L+\omega)^3\\ &= 1-(L+\omega)+(L^2+2L\omega)-3L^2\omega = \big[1-L+L^2\big]+\big[-1+2L-3L^2\big]\omega.\end{aligned}$$

The coefficient in $\omega$ of the integrand of (2.12.1) is therefore

$$(L-KL)(-1+2L-3L^2)+(1-K+2L+c_2-2KL)(1-L+L^2) = 1-K+L^2+c_2;$$

integrated against the fundamental class $[S]$, it gives $L^2+c_2$, as required, see (2.1).

## 2.4 – The Ellingsrud–Göttsche–Lehn recursion

We now describe the inductive computation of [6]. Our aim is to provide the right amount of detail so that: i) one may make sense of the arguments at the end of the proof of Proposition (2.12); and ii) one should be able to make it work in practice. We refer to the original article [6] for the proofs.

Let $S$ be a smooth projective surface, equipped with a line bundle $L \to S$. We consider the Hilbert schemes $S^{[n]}$ and the vector bundles $L^{[n]} \to S^{[n]}$. We call $T_n$ the tangent sheaf of $S^{[n]}$, and $\mathcal{I}_n$ the ideal sheaf of the universal subscheme $Z_n \subseteq S^{[n]} \times S$. For any sequence $I$ of indices in $\{0,\dots,m\}$ we let $\mathrm{pr}_I$ be the projection from $S^{[n]} \times S^m = S^{[n]} \times S \times \cdots \times S$ to the factors indexed by $I$, $S^{[n]}$ being the 0-th factor; beware that we will (abusively) use this notation both for $S^{[n+1]} \times S^m$ and $S^{[n]} \times S^{m+1}$. We let $\Delta$ denote the diagonal in $S^2$.

**(2.14) Proposition.** *Let $f$ be a polynomial in the Chern classes of the following sheaves on $S^{[n+1]} \times S^m$:*

$$\mathrm{pr}_0^*T_{n+1},\ \mathrm{pr}_{0i}^*\mathcal{I}_{n+1},\ \mathrm{pr}_{ij}^*\mathcal{O}_\Delta,\ \mathrm{pr}_i^*T_S,\ \mathrm{pr}_0^*L^{[n+1]},\ \mathrm{pr}_i^*L \qquad (\textit{for } 1 \leqslant i < j \leqslant m).$$

*There is a polynomial $f^\flat$ in the Chern classes of the analogously defined sheaves on $S^{[n]} \times S^{m+1}$, depending only on $f$, and such that*

$$\int_{S^{[n+1]}\times S^m} f = \int_{S^{[n]}\times S^{m+1}} f^\flat.$$

Key to the recursive procedure is the consideration of the incidence variety $S^{[n,n+1]} \subseteq S^{[n]} \times S^{[n+1]}$ of all pairs $(Z, Z')$ such that $Z \subseteq Z'$. There is an identification $S^{[n,n+1]} \cong \mathbf{P}(\mathcal{I}_n)$, coming from the fact that if $Z'$ is obtained by extending $Z$ at the closed point $x \in S$, then the ideal sheaf defining $Z$ in $S$ is an extension of the skyscraper sheaf $\mathbf{C}_x$ by the ideal sheaf of $Z'$.

We let $\mathcal{L} = \mathcal{O}_{\mathbf{P}(\mathcal{I}_n)}(1)$. and consider the following maps: $\phi$ and $\psi$ are the two projections of $\mathbf{P}(\mathcal{I}_n)$ to $S^{[n]}$ and $S^{[n+1]}$, respectively, and $\rho : \mathbf{P}(\mathcal{I}_n) \to S$ maps a pair $(Z, Z')$ to the point $x \in S$ at which $Z$ and $Z'$ differ.

$$\begin{array}{ccccc} & & \mathbf{P}(\mathcal{I}_n) = S^{[n,n+1]} & & \\ & \swarrow^{\phi=\mathrm{pr}_1} & \downarrow{\scriptstyle\rho} & \searrow^{\psi=\mathrm{pr}_2} & \\ S^{[n]} & & S & & S^{[n+1]} \end{array}$$

Moreover, let $j = (\mathrm{id}, \rho) : \mathbf{P}(\mathcal{I}_n) \to \mathbf{P}(\mathcal{I}_n) \times S$, $\phi_S = \phi \times \mathrm{id}_S : \mathbf{P}(\mathcal{I}_n) \times S \to S^{[n]} \times S$, and $\psi_S = \psi \times \mathrm{id}_S : \mathbf{P}(\mathcal{I}_n) \times S \to S^{[n+1]} \times S$. We have the following exact sequence of sheaves on $S^{[n,n+1]} \times S$:

$$0 \to \psi_S^* \mathcal{I}_{n+1} \to \phi_S^* \mathcal{I}_n \to j_* \mathcal{L} \to 0. \tag{2.14.1}$$

Note that $j_*\mathcal{L} = p^*\mathcal{L} \otimes \rho_S^* \mathcal{O}_\Delta$, where $p$ is the projection $S^{[n,n+1]} \times S \to S^{[n,n+1]}$, and $\rho_S = \rho \times \mathrm{id}_S : S^{[n,n+1]} \times S \to S^2$.

*Proof of Proposition (2.14).* We want to relate integrals on $S^{[n+1]} \times S^m$ to integrals on $S^{[n]} \times S^{m+1}$. To this end we consider the product $\Sigma = S^{[n,n+1]} \times S^m$, together with the maps $\Psi = \psi \times \mathrm{id}_{S^m} : \Sigma \to S^{[n+1]} \times S^m$ and $\Phi = (\phi, \rho) \times \mathrm{id}_{S^m} : \Sigma \to S^{[n]} \times S \times S^m$. Moreover, we call $p_\Sigma : \Sigma = S^{[n,n+1]} \times S^m \to S^{[n,n+1]}$ the projection on the first factor.

$$\begin{array}{ccccc} & & \Sigma = S^{[n,n+1]} \times S^m & & \\ & \swarrow^{\Psi=\psi\times\mathrm{id}} & \downarrow{\scriptstyle p_\Sigma} & \searrow^{\Phi=(\phi,\rho)\times\mathrm{id}} & \\ S^{[n+1]} \times S^m & & S^{[n,n+1]} & & S^{[n]} \times S \times S^m \end{array}$$

The first step is to lift our integral to $\Sigma$, using the fact that $\Psi$ is generically finite of degree $n+1$. This gives:

$$\int_{S^{[n+1]}\times S^m} f = \frac{1}{n+1} \int_\Sigma \Psi^* f. \tag{2.14.2}$$

We aim to find an expression for $\Psi^* f$ of the form $\sum_{\nu \geqslant 0} \Phi^* f_\nu \cdot p_\Sigma^*\big(c_1(\mathcal{L})^\nu\big)$, in the Grothendieck ring $K(\Sigma)$ of coherent sheaves on $\Sigma$. Then we will conclude by applying the projection formula.

We first have the following two straightforward formulas (in the Grothendieck ring, as will be the case for all other formulas of the same kind), which display on simple instances the index shift resulting from the additional $S$ factor induced by the $\rho$ component of $\Phi$.

$$\Psi^* \mathrm{pr}_i^* T_S = \Phi^* \mathrm{pr}_{i+1}^* T_S \tag{2.14.3}$$

$$\Psi^* \mathrm{pr}_{ij}^* \mathcal{O}_\Delta = \Phi^* \mathrm{pr}_{i+1,j+1}^* \mathcal{O}_\Delta. \tag{2.14.4}$$

Beware that, by abuse of notation, the 'pr' on the left-hand-side refer to projections from $S^{[n+1]} \times S^m$ whereas those on the right-hand-side refer to projections from $S^{[n]} \times S^{m+1}$. The first non-trivial formula in the induction procedure is the following, and comes from (2.14.1) together with the given expression for $j_*\mathcal{L}$:[1]

$$\Psi^*\mathrm{pr}_{0i}^*\mathcal{I}_{n+1} = \Phi^*\mathrm{pr}_{0,i+1}^*\mathcal{I}_n - p_\Sigma^*\mathcal{L}\cdot\Phi^*\mathrm{pr}_{1,i+1}^*\mathcal{O}_\Delta. \tag{2.14.5}$$

Next we have the following relation, which is a consequence of the computation in the Grothendieck group of $S^{[n,n+1]}$ carried out in [6, §2]:

$$\begin{aligned}(2.14.6)\quad &\Psi^*\mathrm{pr}_0^*T_{n+1}\\ &= \Phi^*\mathrm{pr}_0^*T_n + p_\Sigma^*\mathcal{L}\cdot\Phi^*\mathrm{pr}_{01}^*\mathcal{I}_n^\vee + p_\Sigma^*\mathcal{L}^\vee\cdot\Phi^*\mathrm{pr}_{01}^*\mathcal{I}_n\cdot\Phi^*\mathrm{pr}_1^*\omega_S^\vee - \Phi^*\mathrm{pr}_1^*(\mathcal{O}_S - T_S + \omega_S^\vee).\end{aligned}$$

The last relation follows from [6, Lemma 2.1]:

$$\Psi^*\mathrm{pr}_0^*L^{[n+1]} = \Phi^*\mathrm{pr}_0^*L^{[n]} + p_\Sigma^*\mathcal{L}\cdot\Phi^*\mathrm{pr}_1^*L. \tag{2.14.7}$$

Now, applying formulas (2.14.3), (2.14.4), (2.14.5), (2.14.6), and (2.14.7) to transform $\Psi^*f$, and plugging the result in (2.14.2), one obtains an identity of the form

$$\int_{S^{[n+1]}\times S^m} f = \int_\Sigma\Big(\sum\nolimits_{\nu\geqslant 0}\Phi^*f_\nu\cdot p_\Sigma^*\big(c_1(\mathcal{L})^\nu\big)\Big) \tag{2.14.8}$$

as we wanted. By the projection formula, the right-hand-side equals

$$\int_{S^{[n]}\times S^{m+1}}\Phi_*\Big(\sum\nolimits_{\nu\geqslant 0}\Phi^*f_\nu\cdot p_\Sigma^*\big(c_1(\mathcal{L})^\nu\big)\Big) = \int_{S^{[n]}\times S^{m+1}}\Big(\sum\nolimits_{\nu\geqslant 0}(-1)^\nu f_\nu\cdot c_\nu(-\mathrm{pr}_{01}^*\mathcal{I}_n)\Big),$$

where the last equality is given by the identity $(\phi,\rho)_*\big(c_1(\mathcal{L})^\nu\big) = (-1)^\nu c_\nu(-\mathcal{I}_n)$ in the Chow ring of $S^{[n]}\times S$, see [6, Lemma 1.1].[2] This ends the proof. □

**(2.15) Example.** To illustrate the recursion, let us compute $\int_{S^{[2]}} c_2(L^{[2]})^2$. We consider the diagram

$$\begin{array}{ccccc} & & \Sigma = S^{[1,2]}\times\{\mathrm{pt}\} & & \\ & \swarrow^{\Psi} & \downarrow p & \searrow^{\Phi} & \\ S^{[2]}\times\{\mathrm{pt}\} & & S^{[1,2]} & & S^{[1]}\times S\times\{\mathrm{pt}\}\end{array}$$

in which $\{\mathrm{pt}\}$ denotes the point. This factor may of course be forgotten, except for compatibility with the general notation; thus, the map $p$ is merely the identity. Also, the pair $(S^{[1]}, L^{[1]})$ may be identified with $(S, L)$.

We start from the identity $\int_{S^{[2]}} c_2(L^{[2]})^2 = \frac{1}{2}\int_{S^{[1,2]}}\Psi^*c_2(L^{[2]})^2$. Next, we have $\Psi^*L^{[2]} = \Phi^*\mathrm{pr}_0^*L + p^*\mathcal{L}\cdot\Phi^*\mathrm{pr}_1^*L$ by (2.14.7). In general, $c_2(A+B) = c_2(A)+c_2(B)+c_1(A)c_1(B)$; here, $\Psi^*L^{[2]}$ is expressed as a sum of line bundle, so we simply have

$$c_2(\Psi^*L^{[2]}) = c_1(\Phi^*\mathrm{pr}_0L)\cdot\big(c_1(p^*\mathcal{L}) + c_1(\Phi^*\mathrm{pr}_1^*L)\big)$$

[1] note that there is a typo in the corresponding formula in [6, Proof of Prop. 3.1] (on top of p. 89): the '+' sign there should be a '−' sign; this is the same "small sign error" that is hinted at on the very last line of [18], as was kindly indicated to us by Martijn Kool.

[2] By definition the Chern classes $c_\nu(-\mathcal{F})$ are the coefficients of the multiplicative inverse of the Chern polynomial $c_t(\mathcal{F})$, i.e., the Chern classes $c_\nu(-\mathcal{F})$ are the Segre classes $s_\nu(\mathcal{F})$.

(we have also used that $c_1(AB) = c_1(A) + c_1(B)$), hence

$$c_2(\Psi^* L^{[2]})^2 = \Phi^*\mathrm{pr}_0^* c_1(L)^2 \cdot \big(\underbrace{p^*c_1(\mathcal{L})^2}_{\nu=2} + \underbrace{2p^*c_1(\mathcal{L}) \cdot \Phi^*\mathrm{pr}_1^* c_1(L)}_{\nu=1} + \underbrace{\Phi^*\mathrm{pr}_1^* c_1(L)^2}_{\nu=0}\big)$$

($\nu$ is the index that appears at the end of the proof of Proposition (2.14), from (2.14.8) on).

Now, the universal subscheme $Z_1 \subseteq S^{[1]} \times S$ identifies with the diagonal $\Delta$ in $S \times S$, and we thus have $\mathcal{O}_{S\times S} = \mathcal{I}_1 + \iota_*\mathcal{O}_\Delta$, where $\iota$ denotes the inclusion map $\Delta \hookrightarrow S\times S$. Therefore, $c_1(\mathcal{I}_1) = -c_1(\iota_*\mathcal{O}_\Delta)$, and $c_2(\mathcal{I}_1) = -c_2(\iota_*\mathcal{O}_\Delta)$. By the Riemann–Roch formula without denominators, we have

$$c_1(\iota_*\mathcal{O}_\Delta) = 0 \quad \text{and} \quad c_2(\iota_*\mathcal{O}_\Delta) = -[\Delta],$$

see [7, Example 15.3.1]. Thus,

$$\frac{1}{1 + c_1(\mathcal{I}_1) + c_2(\mathcal{I}_1) + \cdots} = \frac{1}{1 + [\Delta] + \cdots} = 1 - [\Delta] + \cdots ,$$

hence

$$\Phi_* p^* c_1(\mathcal{L}) = -c_1(-\mathcal{I}_1) = 0 \quad \text{and} \quad \Phi_* p^* c_1(\mathcal{L})^2 = c_2(-\mathcal{I}_1) = -[\Delta].$$

In the upshot,

$$\int_{S^{[1,2]}} \Phi^*\mathrm{pr}_0^* c_1(L)^2 \cdot \Phi^*\mathrm{pr}_1^* c_1(L)^2 = \int_S c_1(L)^2 \cdot \int_S c_1(L)^2 = (L^2)^2,$$

$$\int_{S^{[1,2]}} \Phi^*\mathrm{pr}_0^* c_1(L)^2 \cdot 2p^*c_1(\mathcal{L}) \cdot \Phi^*\mathrm{pr}_1^* c_1(L) = -2\int_{S\times S} \mathrm{pr}_0^* c_1(L)^2 \cdot \mathrm{pr}_1^* c_1(L) \cdot c_1(-\mathcal{I}_1) = 0,$$

$$\int_{S^{[1,2]}} \Phi^*\mathrm{pr}_0^* c_1(L)^2 \cdot p^*c_1(\mathcal{L})^2 = \int_{S\times S} \mathrm{pr}_0^* c_1(L)^2 \cdot c_2(-\mathcal{I}_1) = -\int_{S\times S} \mathrm{pr}_0^* c_1(L)^2 \cdot [\Delta] = -L^2,$$

and therefore

$$\int_{S^{[2]}} c_2(L^{[2]})^2 = \frac{1}{2}\int_{S^{[1,2]}} \Psi^* c_2(L^{[2]})^2 = \frac{1}{2}L^2(L^2 - 1).$$

# 3 – Tzeng's proof of the Göttsche–Yau–Zaslow formula

In this section, we outline the proof of Theorem (1.4) given by Tzeng [33]. The plan of the proof is the following. In §3.1, we follow [22] in order to define the algebraic cobordism group $(\omega_{2,1}, +)$ of equivalence classes of pairs $[S, L]$ of surfaces and line bundles with respect to the extended double point relation (cf. Definition (3.3)), endowed with the sum induced by disjoint union of surfaces. Then we recall the description of $\omega_{2,1} \otimes_{\mathbf{Z}} \mathbf{Q}$ as a four dimensional vector space over $\mathbf{Q}$ spanned by $[\mathbf{P}^2, \mathcal{O}_{\mathbf{P}^2}]$, $[\mathbf{P}^2, \mathcal{O}_{\mathbf{P}^2}(1)]$, $[\mathbf{P}^1 \times \mathbf{P}^1, \mathcal{O}_{\mathbf{P}^1\times\mathbf{P}^1}]$, and $[\mathbf{P}^1 \times \mathbf{P}^1, \mathcal{O}_{\mathbf{P}^1\times\mathbf{P}^1}(1,0)]$.

In §3.2 we introduce Göttsche's enumerative integrals $d_\delta(S, L)$ of a line bundle $L$ on a smooth surface $S$, and we show, following Göttsche, that if $L$ is sufficiently ample, then $d_\delta(S, L)$ equals the number $a_\delta(S, L)$ of $\delta$-nodal curves in a general $\delta$-dimensional linear subsystem of $|L|$.

Section 3.3 is dedicated to Tzeng's degeneration formula for the numbers $d_\delta(S, L)$. It can be understood as saying that the map $\phi : (S, L) \mapsto \sum d_\delta(S, L)x^\delta$ is a homomorphism of $\mathbf{Q}$-vector spaces between the algebraic cobordism group $\omega_{2,1} \otimes_{\mathbf{Z}} \mathbf{Q}$ and the space of invertible power series $(\mathbf{Q}[[x]]^\times, \cdot)$.

Finally, §3.4 contains the proofs of Theorems (1.3) and (1.4). The proof of Theorem (1.3) mostly consists in putting together the results of the previous sections; after that, Theorem (1.4) is obtained by plugging in the Bryan–Leung formula in addition.

## 3.1 – The algebraic cobordism group of surfaces and line bundles

We now review the construction of the algebraic cobordism group $(\omega_{2,1},+)$ of line bundles on surfaces, and the results of [22] giving a basis of $\omega_{2,1} \otimes_{\mathbf{Z}} \mathbf{Q}$ as a $\mathbf{Q}$-vector space.

Let $(\mathcal{M}_2,+)$ denote the free abelian group generated by isomorphism classes $[S]$ of smooth projective surfaces. One may think of the sum of two distinct isomorphism classes $[S_1]$ and $[S_2]$ in $\mathcal{M}_2$ as the isomorphism class of the disjoint union $S_1 \sqcup S_2$.

**(3.1) Definition.** *We say that four smooth projective surfaces $S_0,\ldots,S_3$ satisfy a* double point relation $[S_0]-[S_1]-[S_2]+[S_3]$ *if there exists a family $\mathcal{X} \stackrel{\pi}{\longrightarrow} \mathbf{P}^1$ of surfaces such that:*

*(i) the total space $\mathcal{X}$ is smooth of pure dimension 3;*
*(ii) the fibre over 0 is $S_0$;*
*(iii) $S_1$ and $S_2$ are smooth divisors of $\mathcal{X}$ intersecting transversely, and the fibre over $\infty$ is the union $S_1 \cup S_2$;*
*(iv) let $D = S_1 \cap S_2$; the surface $S_3$ is the $\mathbf{P}^1$-bundle $\mathbf{P}(\mathcal{O}_D \oplus N_{D/S_1}) \cong \mathbf{P}(N_{D/S_2} \oplus \mathcal{O}_D)$.*

In the above, $N_{D/S}$ denotes the normal bundle of $D$ in $S$. Note that the surfaces $S_1,S_2$ and the curve $D$ are allowed to be disconnected and empty.

**(3.2) Remark.** The isomorphism $\mathbf{P}(\mathcal{O}_D \oplus N_{D/S_1}) \cong \mathbf{P}(N_{D/S_2} \oplus \mathcal{O}_D)$ in (iv) above follows from the fact that $\mathcal{O}_D(S_1+S_2) \cong \mathcal{O}_D$. Indeed, it implies $N_{D/S_1} \otimes N_{D/S_2} \cong \mathcal{O}_D$, hence $\mathcal{O}_D \oplus N_{D/S_1} \cong (N_{D/S_1} \otimes N_{D/S_2}) \oplus (N_{D/S_1} \otimes \mathcal{O}_D) \cong N_{D/S_1} \otimes (N_{D/S_2} \oplus \mathcal{O}_D)$. See also the discussion of the triple point formula and its variants in [I, Section **].

Analogously, let $(\mathcal{M}_{2,1},+)$ be the free abelian group generated by isomorphism classes $[S,L]$, where $S$ is a smooth projective surfaces and $L$ is a line bundle on $S$, and extend Definition (3.1) as follows.

**(3.3) Definition.** *For $i = 0,\ldots,3$, let $[S_i,L_i]$ be the class of a pair of a smooth projective surface $S_i$ and a line bundle $L_i$ on it. We say that*

$$[S_0,L_0]-[S_1,L_1]-[S_2,L_2]+[S_3,L_3] \tag{3.3.1}$$

*is an* extended double point relation *if $[S_0]-[S_1]-[S_2]+[S_3]$ is a double point relation via a family $\mathcal{X} \stackrel{\pi}{\longrightarrow} \mathbf{P}^1$ of surfaces as above, and there exists a line bundle $\mathcal{L}$ on $\mathcal{X}$ such that $L_i = \mathcal{L}|_{S_i}$ for $i=0,1,2$, and $L_3 = \eta^*(\mathcal{L}|_D)$, where $\eta$ is the structure map*

$$\eta : S_3 = \mathbf{P}(\mathcal{O}_D \oplus N_{D/S_1}) \longrightarrow D.$$

The following fundamental example is closely related to the degeneration to the normal cone. For $\mathbf{P}^2$, it is the degeneration to $\mathbf{P}^2 \cup_{\mathbf{P}^1} \mathbf{F}_1$ used in [3], see [VI].

**(3.4) Example** (see [33, Lemma 2.4])**.** Let us consider a smooth projective surface $S$ and a line bundle $L$ on $S$. Given a smooth curve $C \subseteq S$, let $\mathcal{X} \stackrel{\mathrm{bl}}{\longrightarrow} S \times \mathbf{P}^1$ be the blow-up of the product $S \times \mathbf{P}^1$ along the curve $C \times \{\infty\}$, and let $E \cong \mathbf{P}(\mathcal{O}_C \oplus N_{S/C})$ be its exceptional divisor endowed with the $\mathbf{P}^1$-bundle map $E \stackrel{\eta}{\longrightarrow} C$. The following is an extended double point relation for all integers $k$,

$$\big[S,L\big]-\big[S,L(-kC)\big]-\big[E,\eta^*(L|_C)\otimes\mathcal{O}_E(kC)\big]+\big[E,\eta^*L(-kC)|_C)\big].$$

To see this fact, let $\mathrm{pr}_1 \colon S \times \mathbf{P}^1 \to S$ and $\mathrm{pr}_2 \colon S \times \mathbf{P}^1 \to \mathbf{P}^1$ denote the natural projections. Consider the family $\mathcal{X} \stackrel{\pi}{\longrightarrow} \mathbf{P}^1$ given by $\pi := \mathrm{pr}_2 \circ \mathrm{bl}$, and define a line bundle $\mathcal{L}$ on $\mathcal{X}$ as

$$\mathcal{L} := (\mathrm{bl}^*\mathrm{pr}_1^* L) \otimes \mathcal{O}_{\mathcal{X}}(-kC).$$

The morphism $\pi$ is a smooth over $\mathbf{P}^1 - \{\infty\}$, with fibres $\pi^{-1}(0) \cong S$ and $\pi^{-1}(\infty) \cong S \cup_C E$. Finally, defining $S_0, S_1 \cong S$ and $S_2, S_3 \cong E$ as in Definition (3.1), the line bundle $\mathcal{L}$ is such that $\mathcal{L}|_{S_0} = L$, $\mathcal{L}|_{S_1} = L \otimes \mathcal{O}_S(-kC)$, $\mathcal{L}|_{S_2} = \eta^*(L|_C) \otimes \mathcal{O}_E(kC)$, and $\mathcal{L}|_{S_3} = \eta^*((\mathcal{O}_S(-kC) \otimes L)|_C)$.

We may now define the algebraic cobordism group of line bundles on surfaces.

**(3.5) Definition.** *Let $\mathcal{R}$ denote the subgroup of $\mathcal{M}_{2,1}$ generated by all extended double point relations. The* algebraic cobordism group *of surfaces and line bundles is the quotient*

$$\omega_{2,1} := \mathcal{M}_{2,1}/\mathcal{R}.$$

One may define more generally the algebraic cobordism group $\omega_{n,r}$ of isomorphism classes of pairs of $n$-dimensional varieties and rank $r$-vector bundles. Lee and Pandharipande [22] described the structure of the $\mathbf{Q}$-vector space $\omega_{n,r} \otimes_{\mathbf{Z}} \mathbf{Q}$, by exhibiting a basis and showing that Chern invariants respect algebraic cobordism. In the case of surfaces and line bundles, [22, Theorems 1 and 4] give the following, see [33, Theorem 2.3].

**(3.6) Theorem.** *The map $[S,L] \mapsto (L^2, LK_S, K_S^2, c_2(S))$ is an isomorphism $\omega_{2,1} \otimes_{\mathbf{Z}} \mathbf{Q} \cong \mathbf{Q}^4$ of $\mathbf{Q}$-vector spaces. The family*

$$\big( [\mathbf{P}^2, \mathcal{O}_{\mathbf{P}^2}],\ [\mathbf{P}^2, \mathcal{O}_{\mathbf{P}^2}(1)],\ [\mathbf{P}^1 \times \mathbf{P}^1, \mathcal{O}_{\mathbf{P}^1 \times \mathbf{P}^1}],\ [\mathbf{P}^1 \times \mathbf{P}^1, \mathcal{O}_{\mathbf{P}^1 \times \mathbf{P}^1}(1,0)] \big) \tag{3.6.1}$$

*is a basis of the vector space $\omega_{2,1} \otimes_{\mathbf{Z}} \mathbf{Q}$, and for any pair $[S,L]$, we have*

$$\big[S,L\big] = a_1\big[\mathbf{P}^2, \mathcal{O}_{\mathbf{P}^2}\big] + a_2\big[\mathbf{P}^2, \mathcal{O}_{\mathbf{P}^2}(1)\big] + a_3\big[\mathbf{P}^1 \times \mathbf{P}^1, \mathcal{O}_{\mathbf{P}^1 \times \mathbf{P}^1}\big] + a_4\big[\mathbf{P}^1 \times \mathbf{P}^1, \mathcal{O}_{\mathbf{P}^1 \times \mathbf{P}^1}(1,0)\big],$$

*with $a_1 = -L^2 + \frac{K_S^2 - 2c_2(S)}{3}$, $a_2 = L^2$, $a_3 = \frac{3L^2 + LK_s}{2} - \frac{K_S^2 - 3c_2(S)}{4}$ and $a_4 = -\frac{3L^2 + LK_S}{2}$.*

The first assertion is the content of [22, Theorems 1 and 4], the proof of which we will not discuss. A fundamental input behind it is [20, Theorem 1.2.19].

This being established, proving that the family in (3.6.1) is a basis is elementary: it is enough to check that the 4-tuples $(L^2, LK_S, K_S^2, c_2(S)) \in \mathbf{Q}^4$ corresponding to the basis vectors are linearly independent, which may be seen by computing them explicitly; one finds $(0,0,9,3)$, $(1,-3,9,3)$, $(0,0,8,4)$, and $(0,-2,8,4)$, respectively. Then, one obtains the values of $a_1, \ldots, a_4$ for $[S,L]$ by solving for them in the equality $(L^2, LK_S, K_S^2, c_2(S)) = a_1(0,0,9,3) + \cdots + a_4(0,-2,8,4)$.

**(3.7) Remark.** For future use, note that a family $([S_i, L_i])_{i=1,\ldots,4}$ gives a basis of $\omega_{2,1} \otimes_{\mathbf{Z}} \mathbf{Q}$ if and only if the corresponding vectors $(L_i^2, L_i K_{S_i}, K_{S_i}^2, c_2(S_i)) \in \mathbf{Q}^4$ are linearly independent over $\mathbf{Q}$.

**(3.8) Remark.** Note that the data of the 4-tuples $(L^2, LK_S, K_S^2, c_2(S))$ and $(\chi(L), LK_S, K_S^2, \chi(\mathcal{O}_S))$ are equivalent since, by Noether's formula and the Riemann–Roch theorem,

$$\chi(\mathcal{O}_S) = \frac{1}{12}\left(K_S^2 + c_2(S)\right) \quad \text{and} \quad \chi(L) = \chi(\mathcal{O}_S) + \frac{1}{2}\left(L^2 - LK_S\right).$$

Observe moreover that the above formulas ensure that $L^2 - LK_S$ is even and $K_S^2 + c_2(S)$ is divisible by 3 and 4, so that the coefficients $a_1, \ldots, a_4$ in Theorem (3.6) turn out to be integers, as expected.

**(3.9) Example.** Let $[S_1, L_1]$ and $[S_2, L_2]$ be two classes consisting of a $K3$ surface endowed with a primitive line bundle. Then the 4-tuples $((\chi(L), LK_S, K_S^2, \chi(\mathcal{O}_S))$ corresponding to the family

$$\big( [\mathbf{P}^2, \mathcal{O}_{\mathbf{P}^2}],\ [\mathbf{P}^2, \mathcal{O}_{\mathbf{P}^2}(1)],\ [S_1, L_1],\ [S_2, L_2] \big) \tag{3.9.1}$$

are $(1,0,9,1)$, $(3,-3,9,1)$, $\big(2 + \frac{L_1^2}{2}, 0, 0, 2\big)$ and $\big(2 + \frac{L_2^2}{2}, 0, 0, 2\big)$. This implies that the family in (3.9.1) is a basis of $\omega_{2,1} \otimes_{\mathbf{Z}} \mathbf{Q}$ if and only if $L_1^2 \neq L_2^2$.

### 3.2 – Göttsche's enumerative integrals

In this section we introduce, following Göttsche [8, Section 5], some intersection numbers $d_\delta(S,L)$ on Hilbert schemes of points on $S$ which will be used as a substitute for the numbers $a_\delta(S,L)$ of $\delta$-nodal curves in a general $\delta$-dimensional linear subsystem of $|L|$; this is possible thanks to Theorem (3.11), due to Göttsche as well, which asserts that $a_\delta(S,L)$ and $d_\delta(S,L)$ coincide if $L$ is sufficiently positive.

This substitution circumvents the issue that the numbers $a_\delta(S,L)$ behave rather poorly under degeneration. For instance, we have seen in the first chapters of this volume, e.g., [I] and [VI], that in a degeneration of $S$ to a transverse union of two surfaces $S_1$ and $S_2$, in order to find the number $a_\delta(S,L)$ on $S_1 \cup S_2$ one ought to count not only curves with nodes on $S_1$ and $S_2$, but also curves with tacnodes along the double curve $S_1 \cap S_2$. In contrast, Tzeng's results presented in the next section establish good properties of the numbers $d_\delta(S,L)$.

For any integer $k > 0$, let us consider the Hilbert scheme $S^{[k]}$ parameterizing 0-dimensional subschemes of $S$ of length $k$, and let $Z_k \subseteq S \times S^{[k]}$ be the universal subscheme, with projections $\pi_1\colon Z_k \to S$ and $\pi_2\colon Z_k \to S^{[k]}$. As in §2.3, we denote by $L^{[k]}$ the vector bundle of rank $k$ on $S^{[k]}$ given by $L^{[k]} := \pi_{2*}{\pi_1}^*L$.

**(3.10) Definition** (Göttsche's enumerative integrals). *Let $\mathring{W}^{3\delta} \subseteq S^{[3\delta]}$ be the locally closed subset parametrizing all closed subchemes of the form $\coprod_{i=1}^{\delta} \operatorname{Spec}(\mathcal{O}_{S,x_i}/\mathfrak{m}^2_{S,x_i})$, where $x_1,\dots,x_\delta$ are pairwise distinct points of $S$. Let $W^{3\delta} \subseteq S^{[3\delta]}$ be the closure of $\mathring{W}^{3\delta}$; it is birational to $S^{[\delta]}$. We let*

$$d_\delta(S,L) := \int_{W^{3\delta}} c_{2\delta}(L^{[3\delta]}),$$

*and write the associated generating function*

$$\phi(S,L) := \sum_{\delta \geqslant 0} d_\delta(S,L) x^\delta.$$

In [8, Proposition 5.2], Göttsche proved that if $L$ is sufficiently ample with respect to $\delta$, then the number $d_\delta(S,L)$ coincides with the number $a_\delta(S,L)$ of $\delta$-nodal curves in a general $\delta$-dimensional linear subsystem of $|L|$. Having now Theorem (2.2) at hand, we can slightly weaken the positivity condition and simplify the proof.

**(3.11) Theorem.** *Let $S$ be a smooth projective surface, and let $L$ be a $(3\delta-1)$-very ample line bundle on $S$, for some $\delta > 0$. Then a general linear subsystem of $|L|$ of dimension $\delta$ contains exactly $d_\delta(S,L)$ curves having exactly $\delta$ nodes as singularities, i.e., $a_\delta(S,L) = d_\delta(S,L)$.*

*Proof.* We will prove that a general linear subsystem $\Lambda \subseteq |L|$ of dimension $\delta$ contains $d_\delta(S,L)$ curves (counted with multiplicities) having at least $\delta$ singularities. Then, the result will follow directly from Theorem (2.2), see Remark (2.3).

Let $(s_1,\dots,s_{\delta+1})$ be a basis of the linear subspace of $H^0(S,L)$ corresponding to $\Lambda \subseteq \mathbf{P}(H^0(S,L))$, and consider the global sections $\pi_{2*}\pi_1^*s_1,\dots,\pi_{2*}\pi_1^*s_{\delta+1}$ of the vector bundle $L^{[3\delta]}$. Moreover, let $Y \subseteq S^{[3\delta]}$ be the locus of points $[Z] \in S^{[3\delta]}$ such that the vectors $\pi_{2*}{\pi_1}^*s_1(Z),\dots,\pi_{2*}{\pi_1}^*s_{\delta+1}(Z)$ are linearly dependent, i.e., such that there exists a divisor $D$ belonging to $\Lambda$ that contains $Z$.

Since $L$ is $(3\delta-1)$-very ample, the restriction map $H^0(S,L) \to H^0(Z, L\otimes\mathcal{O}_Z)$ is surjective for all $[Z] \in S^{[3\delta]}$, hence the evaluation map $H^0(S,L)\otimes\mathcal{O}_{S^{[3\delta]}} \to L^{[3\delta]}$ is surjective as well, and $L^{[3\delta]}$ is globally generated. Thus, by [7, Example 14.3.2] applied to the restriction of $L^{[3\delta]}$ to $W^{3\delta}$, we find that $Y \cap W^{3\delta}$ has codimension $2\delta$ in $W^{3\delta}$, and class $[Y\cap W^{3\delta}] = c_{2\delta}(L^{[3\delta]}) \cap [W^{3\delta}]$ by the Thom–Porteous formula. Moreover, since the dimension of $W^{3\delta} \setminus \mathring{W}^{3\delta}$ is smaller than $2\delta$, the same argument applied to the restriction of $L^{[3\delta]}$ to $W^{3\delta} \setminus \mathring{W}^{3\delta}$ gives that $Y \cap W^{3\delta}$ lies entirely in $\mathring{W}^{3\delta}$.

We conclude that the number of members of $\Lambda$ (counted with multiplicities) containing a point of $\mathring{W}^{3\delta}$ is $d_\delta(S,L) := \int_{W^{3\delta}} c_{2\delta}(L^{[3\delta]})$. Now a curve $C \subseteq S$ contains a scheme $Z \in \mathring{W}^{3\delta}$ with $Z_{\mathrm{red}} = x_1 + \cdots + x_\delta$ if and only if it is singular at all points $x_1, \ldots, x_\delta \in S$. Therefore, the number of curves in $\Lambda$ having at least $\delta$ singularities is $d_\delta(S,L)$, and this ends the proof. □

## 3.3 – Degeneration of enumerative integrals

We now follow [33, Section 3] in order to show how the generating function $\phi(S,L) = \sum d_\delta(S,L)x^\delta$ of the Göttsche enumerative integrals of a pair $(S,L)$ respects extended double point relations. The result below is Tzeng's main contribution to the proof of the Göttsche conjecture. It rests on Proposition (3.15), which may be seen as a version of Jun Li's celebrated degeneration formula [23].

**(3.12) Theorem.** *If* $[S_0,L_0] - [S_1,L_1] - [S_2,L_2] + [S_3,L_3]$ *is an extended double point relation in* $\mathcal{M}_{2,1}$, *then*

$$\phi(S_0,L_0) = \frac{\phi(S_1,L_1)\phi(S_2,L_2)}{\phi(S_3,L_3)}.$$

*Equivalently, the map* $\phi\colon \big(\omega_{2,1} \otimes_{\mathbf{Z}} \mathbf{Q}, +\big) \longrightarrow \big(\mathbf{Q}[[x]]^\times, \times\big)$ *is a* $\mathbf{Q}$*-linear group homomorphism.*

Let us consider a family $\mathcal{X} \xrightarrow{\pi} \mathbf{P}^1$ of surfaces and a line bundle $\mathcal{L}$ inducing the extended double point relation $[S_0,L_0] - [S_1,L_1] - [S_2,L_2] + [S_3,L_3]$ appearing in the statement, in the same notation as in Definitions (3.1) and (3.3). For all $t \in \mathbf{P}^1$, we let $X_t$ be the fibre $\pi^{-1}(t)$. The fibre $X_0$ is $S_0$, and the fibre $X_\infty$ is the union of $S_1$ and $S_2$ meeting transversally along a smooth curve $D$.

Let $\delta$ be a positive integer; for all smooth surfaces $S$ we let $W_S^{3\delta} \subseteq S^{[3\delta]}$ be as in Definition (3.10). The strategy of [33] to prove Theorem (3.12) is to study the degeneration of $W_{X_t}^{3\delta}$ at $t = \infty$ by means of the good degeneration of Hilbert schemes constructed by Li and Wu [25], see also [24, Section 3.4].

Let $U \subseteq \mathbf{P}^1$ be the largest open subset containing $\infty$ and such that $\pi$ is smooth over $U - \{\infty\}$, and let us still denote by $\mathcal{X} \xrightarrow{\pi} U$ the pullback of the original family over $\mathbf{P}^1$. For all $n \geqslant 1$, Li and Wu construct a family $\pi^{[n]}\colon \mathcal{X}^{[n]} \to U$ which is a proper and separated Deligne–Mumford stack of finite type over $U$, such that for all $t \neq \infty$, the fibre of $\mathcal{X}^{[n]}$ over $t$ is the Hilbert scheme $X_t^{[n]}$ of $n$ points on the smooth surface $X_t$; the fibre over $\infty$ is different, however, as we will see.

The first step in the construction is to replace $\mathcal{X} \xrightarrow{\pi} U$ by an alternative family $\pi[n]\colon \mathcal{X}[n] \to U$ of surfaces, called an *expansion*, equal to $\mathcal{X} \to U$ over $U - \{\infty\}$, and with fibre over $\infty$ a semistable model of $S_1 \cup S_2$, viz.

$$\mathcal{X}[n]_\infty := S_1 \cup \Delta_1 \cup \Delta_2 \cup \cdots \cup \Delta_{n-1} \cup S_2$$

where, for all $i = 1, \ldots, n-1$, the surface $\Delta_i$ is a $\mathbf{P}^1$-bundle $\mathbf{P}(\mathcal{O}_D \oplus N_{S_1/D})$ over $D$ with two distinguished sections $D_i^-$ and $D_i^+$, and $\mathcal{X}[n]_\infty$ is obtained by identifying $D_i^+ \subseteq \Delta_i$ with $D_{i+1}^- \subseteq \Delta_{i+1}$, $D \subseteq S_1$ with $D_1^- \subseteq \Delta_1$, and $D_{n-1}^+ \subseteq \Delta_{n-1}$ with $D \subseteq S_2$.

Li and Wu introduce the notion of *stable ideal sheaves*, see [24, Definition 3.53], and let $\pi^{[n]}\colon \mathcal{X}^{[n]} \to U$ be the moduli space of stable ideal sheaves of the family $\mathcal{X}[n] \to U$. For all $t \neq \infty$, the fibre ${\pi^{[n]}}^{-1}(t)$ is the Hilbert scheme $X_t^{[n]}$ parametrizing 0-dimensional subschemes of $X_t$ of length $n$, whereas the fibre ${\pi^{[n]}}^{-1}(\infty)$ is the union

$$\mathcal{X}_\infty^{[n]} := \bigcup_{k=0}^{n} \big((S_1/D)^{[k]} \times (S_2/D)^{[n-k]}\big), \tag{3.12.1}$$

where $(S_1/D)^{[k]}$ (resp. $(S_2/D)^{[k]}$) parametrizes 0-dimensional schemes of length $k$ supported on the smooth locus of $S_1 \cup \Delta_1 \cup \cdots \cup \Delta_i$ for some $i \leqslant k$ (resp. $\Delta_j \cup \cdots \cup \Delta_{n-1} \cup S_2$ for some $j \leqslant n-k$);[3] by convention, $(S_1/D)^{[0]} = (S_2/D)^{[0]} = \{\mathrm{pt}\}$.

We now set $n = 3\delta$, and consider the expanded family of surfaces $\pi[3\delta] \colon \mathcal{X}[3\delta] \to U$, and the moduli space of stable ideal sheaves $\pi^{[3\delta]} \colon \mathcal{X}^{[3\delta]} \to U$.[4] The following lemma describes the good degeneration of $W^{3\delta}_{X_t}$ in $\mathcal{X}^{[3\delta]}$ at $t = \infty$. Let $\mathcal{W}^{3\delta}$ be the the closure in $\mathcal{X}^{[3\delta]}$ of the set

$$\mathring{\mathcal{W}}^{3\delta} := \coprod_{t \neq \infty} \mathring{W}^{3\delta}_{X_t} = \Big\{ \Big[ \coprod_{i=1}^{\delta} \mathrm{Spec}\left(\mathcal{O}_{X_t,x_i}/\mathfrak{m}^2_{x_i}\right) \Big] : t \neq \infty \text{ and } x_1, \ldots, x_\delta \in X_t \text{ pairwise distinct} \Big\}.$$

For all $k = 1, \ldots, \delta$ and $i = 1, 2$, we let $W^{3k}_{S_i/D}$ be the closure in $(S_i/D)^{[3k]}$ of the locus parametrizing unions of $k$ disjoint subschemes of the form $\mathrm{Spec}(\mathcal{O}_x/\mathfrak{m}^2_x)$; by convention, $W^0_{S_i/D} := \{\mathrm{pt}\}$.

**(3.13) Lemma** ([33, Lemma 3.8])**.** *The map $\mathcal{W}^{3\delta} \to U$ is flat, its fibre over $t \neq \infty$ is $W^{3\delta}_{X_t}$, and its fibre over $\infty$ is such that for all $m = 0, \ldots, 3\delta$,*

$$\mathcal{W}^{3\delta} \cap \left( (S_1/D)^{[m]} \times (S_2/D)^{[3\delta - m]} \right) = \begin{cases} W^{3k}_{S_1/D} \times W^{3(\delta-k)}_{S_2/D} & \text{if } m = 3k \text{ for some } k \in \mathbf{N} \\ \varnothing & \text{otherwise.} \end{cases}$$

Let us now briefly survey the definitions of a universal subscheme over $\mathcal{X}^{[n]}$, and of a vector bundle $\mathcal{L}^{[n]}$ on $\mathcal{X}^{[n]}$ from the line bundle $\mathcal{L}$ on $\mathcal{X}$, see [33, Definition-Proposition 3.4]. Let $n$ be a positive integer. We consider the universal closed subscheme $Z^{[n]}_0 \subseteq S_0 \times S^{[n]}_0$. For $i = 1, 2$ and $0 \leqslant k \leqslant n$, we may also introduce universal closed subschemes $Z^{[k]}_i \subseteq (S_i/D) \times (S_i/D)^{[k]}$; beware that $S_i/D$ denotes a suitable expansion of $S_i$ along the curve $D$. Let

$$\begin{array}{ccccccc} Z^{[n]}_0 & \xrightarrow{q^{[n]}_0} & S^{[n]}_0 & \quad\text{and}\quad & Z^{[k]}_i & \xrightarrow{q^{[k]}_i} & (S_i/D)^{[k]} \\ {\scriptstyle p^{[n]}_0}\downarrow & & & & {\scriptstyle p^{[k]}_i}\downarrow & & \\ S_0 & & & & (S_i/D) & & \end{array}$$

be the corresponding projections, and define the vector bundles of rank $3k$

$$L^{[k]}_i := (q^{[k]}_i)_* (p^{[k]}_i)^* L.$$

for all $i = 0, 1, 2$. We denote by $\pi^{[k,n]}_1$ and $\pi^{[k,n]}_2$ the natural projections

$$\begin{array}{ccc} (S_1/D)^{[k]} \times (S_2/D)^{[n-k]} & \xrightarrow{\pi^{[k,n]}_2} & (S_2/D)^{[n-k]} \\ {\scriptstyle \pi^{[k,n]}_1}\downarrow & & \\ (S_1/D)^{[k]}. & & \end{array}$$

[3]For a subscheme $Z$ supported on $S_1 \cup \Delta_1 \cup \cdots \cup \Delta_i$ (resp. $\Delta_j \cup \cdots \cup \Delta_{n-1} \cup S_2$), we require moreover that all components $\Delta_1, \ldots, \Delta_i$ (resp. $\Delta_j, \ldots, \Delta_{n-1}$) contain a point of the support of $Z$ (the components not supporting any piece of $Z$ ought to be contracted).

[4]beware the unfortunate proximity of the two otherwise standard notations.

With this notation, there exists a universal closed subscheme $\mathcal{Z} \subseteq \mathcal{X}[n] \times \mathcal{X}^{[n]}$ with projections

$$\begin{array}{ccc} \mathcal{Z} & \xrightarrow{Q} & \mathcal{X}^{[n]} \\ {\scriptstyle P}\downarrow & & \\ \mathcal{X}[n]. & & \end{array}$$

Its fibres over 0 and $\infty$ are $Z_0^{[n]}$ and

$$\mathcal{Z}_\infty := \left( \bigcup_{k=0}^{n} \left( Z_1^{[k]} \times (S_2/D)^{[n-k]} \right) \right) \cup \left( \bigcup_{k=0}^{n} \left( (S_1/D)^{[k]} \times Z_2^{[n-k]} \right) \right),$$

respectively.

There is a line bundle $\mathcal{L}[n]$ on $\mathcal{X}[n]$ obtained by pull-back from the line bundle $\mathcal{L}$ on $\mathcal{X}$. Setting $\mathcal{L}^{[n]} := Q_* P^* \mathcal{L}[n]$, the following holds.

**(3.14) Lemma** ([33, Lemmas 3.6 and 3.7]). *The sheaf $\mathcal{L}^{[n]}$ is a vector bundle of rank $3n$ over $\mathcal{X}^{[n]}$, such that*

$$\mathcal{L}^{[n]}|_{S_0} = L_0^{[n]} \quad \textit{and} \quad \mathcal{L}^{[n]}|_{(S_1/D)^{[k]} \times (S_2/D)^{[n-k]}} = (\pi_1^{[k,n]})^* L_1^{[k]} \oplus (\pi_2^{[k,n]})^* L_2^{[n-k]}.$$

In order to study degenerations of the Göttsche enumerative numbers, we define for $i = 1, 2$ the *relative Göttsche integrals*

$$d_k(S_i/D, L_i) = \begin{cases} 1 & \text{for } k = 0 \\ \displaystyle\int_{W^{3k}_{S_i/D}} c_{2k}(L_i^{[3k]}) & \text{for } k \geqslant 1, \end{cases}$$

and their *relative generating function*

$$\phi(S_i/D, L_i) = \sum_{k \geqslant 0} d_k(S_i/D, L_i) x^k.$$

Note that the relative Göttsche integrals are defined independently of the family $\mathcal{X} \to U$, and depend only on $S_i$ and $D \subseteq S_i$: to define $(S_i/D)^{[3k]}$ and $L_i^{[3k]}$, one considers the degeneration to the normal cone as in Example (3.4).

We may now finally state the degeneration formula which lies at the heart of Tzeng's proof of the Göttsche conjecture.

**(3.15) Proposition** ([33, Proposition 3.9]). *If $[S_0, L_0] - [S_1, L_1] - [S_2, L_2] + [S_3, L_3]$ is an extended double point relation in $\mathcal{M}_{2,1}$, then*

(3.15.1)
$$d_\delta\left(S_0, L_0\right) = \sum_{k=0}^{\delta} d_k\left(S_1/D, L_1\right) d_{\delta-k}\left(S_2/D, L_2\right),$$

*hence*

(3.15.2)
$$\phi(S_0, L_0) = \phi(S_1/D, L_1) \cdot \phi(S_2/D, L_2).$$

We emphasize that while $d_\delta(S_0, L_0)$ is a genuine Göttsche integral, as defined in Section 3.2, the $d_k(S_i/D, L_i)$ are relative, i.e., they are defined in terms of good degenerations à la Jun Li et al.

*Proof.* Formula (3.15.1) follows from the good degeneration results surveyed above, after fairly mechanical computations. The moduli space of stable ideal sheaves $\mathcal{X}^{[3\delta]} \to U$ provides a degeneration of the Hilbert scheme $S_0^{[3\delta]}$ to a moduli space defined on an expansion of $S_1 \cup_D S_2$, see (3.12.1). Within this degeneration, the family $\mathcal{W}^{3\delta} \to U$ provides a flat degeneration of $W^{3\delta}_{S_0}$ to the fibre $\mathcal{W}^{3\delta}_\infty$ over $\infty \in U$, see Lemma (3.13), and the vector bundle $\mathcal{L}^{[3\delta]}$ provides a degeneration of $L_0^{[3\delta]}$, see Lemma (3.14). Formula (3.15.1) will follow from the equality

$$c_{2\delta}(\mathcal{L}^{[3\delta]})\big[\mathcal{W}_0^{[3\delta]}\big] = c_{2\delta}(\mathcal{L}^{[3\delta]})\big[\mathcal{W}_\infty^{[3\delta]}\big] :$$

the left-hand-side plainly writes

$$c_{2\delta}(\mathcal{L}^{[3\delta]})\big[\mathcal{W}_0^{[3\delta]}\big] = c_{2\delta}(\mathcal{L}^{[3\delta]})\big[W_{S_0}^{[3\delta]}\big] = \int_{W^{3\delta}_{S_0}} c_{2\delta}(L_0^{[3\delta]}) = d_\delta(S_0, L_0),$$

and we will see that the right-hand-side equals the right-hand-side of (3.15.1), which will end the proof.

By Lemmas (3.13) and (3.14), we have

$$\begin{aligned} c_{2\delta}(\mathcal{L}^{[3\delta]})\big[\mathcal{W}_\infty^{[3\delta]}\big] &= c_{2\delta}(\mathcal{L}^{[3\delta]}) \Big[\bigcup\nolimits_{k=0}^{\delta} W^{3k}_{S_1/D} \times W^{3(\delta-k)}_{S_2/D}\Big] \\ &= \sum_{k=0}^{\delta} c_{2\delta}\big( (\pi_1^{[3k,3\delta]})^* L_1^{[3k]} \oplus (\pi_2^{[3k,3\delta]})^* L_2^{[3(\delta-k)]} \big) \big[W^{3k}_{S_1/D} \times W^{3(\delta-k)}_{S_2/D}\big]. \end{aligned}$$

Now,

$$c_{2\delta}\big((\pi_1^{[3k,3\delta]})^* L_1^{[3k]} \oplus (\pi_2^{[3k,3\delta]})^* L_2^{[3(\delta-k)]}\big) = \sum_{l=0}^{2\delta} c_l\big((\pi_1^{[3k,3\delta]})^* L_1^{[3k]}\big) c_{2\delta-l}\big((\pi_2^{[3k,3\delta]})^* L_2^{[3(\delta-k)]}\big),$$

so that after capping with $\big[W^{3k}_{S_1/D} \times W^{3(\delta-k)}_{S_2/D}\big]$ the only non-zero term is obtained for $l = \dim W^{3k}_{S_1/D} = 2k$. Therefore, we finally obtain

$$\begin{aligned} c_{2\delta}(\mathcal{L}^{[3\delta]})\big[\mathcal{W}_\infty^{[3\delta]}\big] &= \sum_{k=0}^{\delta} c_{2k}(L_1^{[3k]})\big[W^{3k}_{S_1/D}\big] c_{2(\delta-k)}(L_2^{[3(\delta-k)]})\big[W^{3(\delta-k)}_{S_2/D}\big] \\ &= \sum_{k=0}^{\delta} d_k\left(S_1/D, L_1\right) d_{\delta-k}\left(S_2/D, L_2\right), \end{aligned}$$

which ends the proof of Formula (3.15.1). Formula (3.15.2) follows at once. □

With Proposition (3.15) at hand, the proof of Theorem (3.12) consists in using degenerations to the normal cone (and Proposition (3.15) itself) in order to replace all relative generating functions by absolute ones in the degeneration formula (3.15.2). We will thus use several times the following lemma.

**(3.16) Lemma.** *Let $S$ be a smooth surface, and $C \subseteq S$ be a smooth curve with normal bundle $N_{C/S}$. Let $E \cong \mathbf{P}(\mathcal{O}_C \oplus N_{C/S})$, and $\eta : E \to C$ be the $\mathbf{P}^1$-bundle map. For all line bundles $L$ on $S$, we have*

(3.16.1) $$\phi(S, L) = \phi(S/C, L) \cdot \phi(E/C, \eta^*(L|_C)).$$

*Proof.* Let us consider the degeneration of $S$ to $S \cup_C E$ constructed in Example (3.4), endowed with the line bundle $\mathcal{L}$ obtained by pulling $L$ back (this is the case $k = 0$ in Example (3.4)). On the degenerate fibre $S \cup_C E$, one has $\mathcal{L}|_S = L$ and $\mathcal{L}|_E = \eta^*\, L|_C$. Therefore, the degeneration formula (3.15.2) applied to this situation gives the required Formula (3.16.1). □

**(3.17)** *Proof of Theorem (3.12).* We have

$$\phi(S_0, L_0) = \phi(S_1/D, L_1) \cdot \phi(S_2/D, L_2). \tag{3.17.1}$$

by Proposition (3.15).

By Lemma (3.16), we have

$$\phi(S_i, L_i) = \phi(S_i/D, L_i) \cdot \phi(E_i/D_i^-, \eta_i^* \, L_i|_D) \tag{3.17.2}$$

for $i = 1, 2$, where $E_i := \mathbf{P}(\mathcal{O}_D \oplus N_{D/S_i})$, $\eta_i : E_i \to D$ is the $\mathbf{P}^1$-bundle map, and $D_i^- := \mathbf{P}(N_{X_i/D})$ and $D_i^+ := \mathbf{P}(\mathcal{O}_D)$ are the two distinguished sections of $\eta_i$. Now, the assumed extended double point relation tells us that $S_3 \simeq E_i$, and $L_3 \simeq \eta_i^* \, \mathcal{L}|_D$ for $i = 1, 2$, in the notation introduced after Theorem (3.12). Since $\mathcal{L}|_{S_i} = L_i$, we have $L_3 \simeq \eta_i^* \, L_i|_D$, so that (3.17.2) may be rewritten as

$$\phi(S_i, L_i) = \phi(S_i/D, L_i) \cdot \phi(E_i/D_i^-, L_3). \tag{3.17.3}$$

Finally, we apply Lemma (3.16) to the curve $D_1^-$ in the surface $E_1$; in this case, the degeneration to the normal cone considered in the proof of the lemma is a degeneration of $E_1$ to $E_1 \cup_D E_1$, where the curve $D$ is embedded as $D_1^-$ in the first $E_1$, and as $D_1^+$ in the second $E_1$. We obtain the identity

$$\phi(E_1, L_3) = \phi(E_1/D_1^-, L_3) \cdot \phi(E_1/D_1^+, L_3). \tag{3.17.4}$$

Since $E_1 \simeq S_3$ on the one hand, and on the other hand the isomorphism $(E_1, \eta_1^* \, L_1|_D)) \simeq (E_2, \eta_2^* \, L_2|_D)$ maps $D_1^+$ to $D_2^-$ (see Remark (3.2)), (3.17.4) may be rewritten as

$$\phi(S_3, L_3) = \phi(E_1/D_1^-, L_3) \cdot \phi(E_2/D_2^-, L_3). \tag{3.17.5}$$

In conclusion, by (3.17.1), (3.17.3), and (3.17.5), we have

$$\phi(S_0, L_0) = \frac{\phi(S_1, L_1)}{\phi(E_1/D_1^-, \eta_1^* \, L_1|_D)} \cdot \frac{\phi(S_2, L_2)}{\phi(E_2/D_2^-, \eta_2^* \, L_2|_D)} = \frac{\phi(S_1, L_1) \cdot \phi(S_2, L_2)}{\phi(S_3, L_3)},$$

as we wanted. □

## 3.4 – Conclusion

This final section is dedicated to the proofs of Theorems (1.3) and (1.4). The substantial work has already been done, and it remains only to assemble the various pieces.

We first prove the analogue of Theorem (1.3) for the generating function $\phi(S, L) = \sum d_\delta(S, L) x^\delta$, and describe explicitly the universal power series involved. Theorem (1.3) itself will follow rather directly, thanks to Theorem (3.11).

**(3.18) Proposition.** *For all pairs* $(S, L)$*, one has*

$$\phi(S, L) = A_1^{L^2} A_2^{LK_S} A_3^{K_S^2} A_4^{c_2(S)},$$

*with*

$$A_1 = \phi(\mathbf{P}^2, \mathcal{O}_{\mathbf{P}^2})^{-1} \phi(\mathbf{P}^2, \mathcal{O}_{\mathbf{P}^2}(1)) \phi(\mathbf{P}^1 \times \mathbf{P}^1, \mathcal{O}_{\mathbf{P}^1 \times \mathbf{P}^1})^{\frac{3}{2}} \phi(\mathbf{P}^1 \times \mathbf{P}^1, \mathcal{O}_{\mathbf{P}^1 \times \mathbf{P}^1}(1,0))^{-\frac{3}{2}},$$
$$A_2 = \phi(\mathbf{P}^1 \times \mathbf{P}^1, \mathcal{O}_{\mathbf{P}^1 \times \mathbf{P}^1})^{\frac{1}{2}} \phi(\mathbf{P}^1 \times \mathbf{P}^1, \mathcal{O}_{\mathbf{P}^1 \times \mathbf{P}^1}(1,0))^{-\frac{1}{2}},$$
$$A_3 = \phi(\mathbf{P}^2, \mathcal{O}_{\mathbf{P}^2})^{-\frac{1}{3}} \phi(\mathbf{P}^1 \times \mathbf{P}^1, \mathcal{O}_{\mathbf{P}^1 \times \mathbf{P}^1})^{-\frac{1}{4}},$$
$$A_4 = \phi(\mathbf{P}^2, \mathcal{O}_{\mathbf{P}^2})^{-\frac{2}{3}} \phi(\mathbf{P}^1 \times \mathbf{P}^1, \mathcal{O}_{\mathbf{P}^1 \times \mathbf{P}^1})^{\frac{3}{4}}.$$

*Proof.* Given a class $[S,L] \in \omega_{2,1}$, Theorem (3.6) tells us that

$$[S,L] = a_1[\mathbf{P}^2, \mathcal{O}_{\mathbf{P}^2}] + a_2[\mathbf{P}^2, \mathcal{O}_{\mathbf{P}^2}(1)] + a_3[\mathbf{P}^1 \times \mathbf{P}^1, \mathcal{O}_{\mathbf{P}^1\times\mathbf{P}^1}] + a_4[\mathbf{P}^1 \times \mathbf{P}^1, \mathcal{O}_{\mathbf{P}^1\times\mathbf{P}^1}(1,0)],$$

with $a_1 = -L^2 + \frac{K_S^2 - 2c_2(S)}{3}$, $a_2 = L^2$, $a_3 = \frac{3L^2+LK_s}{2} - \frac{K_S^2-3c_2(S)}{4}$ and $a_4 = -\frac{3L^2+LK_S}{2}$. Then, Theorem (3.12) yields

$$\phi(S,L) = \phi(\mathbf{P}^2, \mathcal{O}_{\mathbf{P}^2})^{a_1} \phi(\mathbf{P}^2, \mathcal{O}_{\mathbf{P}^2}(1))^{a_2} \phi(\mathbf{P}^1 \times \mathbf{P}^1, \mathcal{O}_{\mathbf{P}^1\times\mathbf{P}^1})^{a_3} \phi(\mathbf{P}^1 \times \mathbf{P}^1, \mathcal{O}_{\mathbf{P}^1\times\mathbf{P}^1}(1,0))^{a_4},$$

from which the assertion follows by replacing $a_1, \dots, a_4$ with their expressions in terms of $L^2$, $LK_S$, $K_S^2$, and $c_2(S)$. □

**(3.19) Remark.** According to Remark (3.7), we may easily determine an alternative basis of $\omega_{2,1} \otimes_{\mathbf{Z}} \mathbf{Q}$, and the corresponding coefficients of an arbitrary class $[S,L] \in \omega_{2,1}$ may be otained as in the proof of Theorem (3.6). Consequently, the universal power series $A_1, A_2, A_3, A_4 \in \mathbf{Q}[[x]]$ could be described in terms of any such basis.

**(3.20) Remark.** Tzeng proved the Göttsche conjecture as follows (see [33, Theorem 1.1]). She observed that by expanding the series $A_1, A_2, A_3, A_4$ in Proposition (3.18), the coefficient $d_\delta(S,L)$ of $x^\delta$ in the generating function $\phi(S,L)$ is a universal polynomial $T_\delta$ of degree $\delta$ in the variables $L^2, LK_S, K_S^2, c_2(S)$ (cf. [33, Corollary 4.2]), i.e., the analogue of the Göttsche conjecture holds for the numbers $d_\delta(S,L)$, without any positivity assumption on $L$. Then, the conjecture for the numbers $a_\delta(S,L)$ follows by Theorem (3.11), under a suitable positivity asumption (Tzeng needed $(5\delta-1)$-very ampleness, but we now know that $(3\delta-1)$-very ampleness suffices, thanks to [18], see Theorem (2.2)).

**(3.21)** *Proof of Theorem (1.3).* Consider the generating function $T(S,L)$ of the universal polynomials $T_\delta(L^2, LK_S, K_S^2, c_2(S))$; the existence of the latter is guaranteed by Theorem (1.2). If $L$ is $(3\delta-1)$-very ample, then Theorems (1.2) and (3.11) ensure that $T_\delta(L^2, LK_S, K_S^2, c_2(S)) = a_\delta(S,L) = d_\delta(S,L)$. Thus the generating functions $T(S,L)$ and $\phi(S,L)$ coincide, and we deduce from Proposition (3.18) that $T(S,L) = A_1^{L^2} A_2^{LK_S} A_3^{K_S^2} A_4^{c_2(S)}$. □

Finally, let us prove the Göttsche–Yau–Zaslow formula.

**(3.22)** *Proof of Theorem (1.4).* Let us consider the generating function

$$\gamma(S,L)(q) := \sum_{\delta \geqslant 0} T_\delta(L^2, LK_S, K_S^2, c_2(S))\, (DG_2(\tau))^\delta$$

appearing in the right-hand side of Göttsche–Yau–Zaslow formula. Given two integers $g \geqslant \delta \geqslant 0$, let $(S,L)$ consist of a very general $K3$ surface $S$ endowed with a primitive line bundle $L$ such that $L^2 = 2(g+\delta-1)$, and let $M_g(\delta)$ denote the number of curves in $|L|$ having geometric genus $g$ and passing through $g$ general points of $S$ (or, equivalently, contained in a general linear subsystem of $V \subseteq |L|$ of codimension $g$ and dimension $\delta$). Bryan and Leung [2]) (see [XI, Theorem **]) proved that the generating function of the numbers $M_g(\delta)$ satisfies the formula

$$\sum_{\delta \geqslant 0} M_g(\delta) q^{g+\delta-1} = \frac{DG_2(\tau)^g}{\Delta(\tau)}$$

which, according to [8, Remark 2.6], can be reformulated as

$$\gamma(S,L)(q) = \frac{(DG_2(\tau)/q)^{\chi(L)}}{\Delta(\tau) D^2 G_2(\tau)/q^2}. \tag{3.22.1}$$

On the other hand, we recall that the 4-tuples $(L^2, LK_S, K_S^2, c_2(S))$ and $(K_S^2, LK_S, \chi(L), \chi(\mathcal{O}_S))$ determine each other through linear relations (cf. Remark (3.8)). Since $\gamma(S,L)(q)$ coincides with the generating function $\phi(S,L)$ computed at $x = DG_2(\tau)$, we deduce from Proposition (3.18) that there exist four universal power series $B_1(q), \dots, B_4(q) \in \mathbf{Q}[[q]]$ such that

$$\gamma(S,L)(q) = B_1(q)^{K_S^2} B_2(q)^{LK_S} B_3(q)^{\chi(L)} B_4(q)^{\chi(\mathcal{O}_S)}. \tag{3.22.2}$$

When $S$ is a $K3$ surface, equation (3.22.2) becomes $\gamma(S,L)(q) = B_3(q)^{\chi(L)} B_4(q)^2$. Moreover, as we consider different pairs $(S,L)$ such that $L^2 = 2(g+\delta-1)$, the latter equality must be satisfied for any positive integer $\chi(L) = g+\delta+1$. Comparing $\gamma(S,L)(q) = B_3(q)^{\chi(L)} B_4(q)^2$ and (3.22.1), we then conclude that $B_3(q) = DG_2(\tau)/q$ and $B_4(q) = \big(\Delta(\tau) D^2 G_2(\tau)/q^2\big)^{-1/2}$. Therefore (3.22.2) leads to the Göttsche–Yau–Zaslow formula

$$\gamma(S,L)(q) = B_1(q)^{K_S^2} B_2(q)^{LK_S} \frac{(DG_2(\tau)/q)^{\chi(L)}}{(\Delta(\tau) D^2 G_2(\tau)/q^2)^{\chi(\mathcal{O}_S)/2}}.$$

□

## References

### Chapters in this Volume[a] [b]

### External References

[a] see `https://www.math.univ-toulouse.fr/~tdedieu/#enumTV`

[b] For LaTeX reasons, references to labels in other chapters of the volume have not been calculated in this version, causing references of the form [I, Thm. **] instead of [I, Thm. (1.1)]; we apologize for the inconvenience, and refer to the full volume for appropriately expanded references.